\documentclass[11pt,reqno]{amsart}

\usepackage{amsmath}
 \usepackage{amssymb}

\newcommand{\mysection}[1]{\section{#1}
      \setcounter{equation}{0}}

\newtheorem{theorem}{Theorem}[section]
\newtheorem{lemma}[theorem]{Lemma}

\theoremstyle{definition}
\newtheorem{assumption}{Assumption}[section]
\newtheorem{definition}{Definition}[section]
\theoremstyle{remark}
\newtheorem{remark}{Remark}[section]

\newtheorem{example}{Example}[section]

\renewcommand\){{\rm)}}

\newcommand\bbeta{\text{\raise-.2ex\hbox{$\bm{\beta}$}}}

\makeatletter
 \def\dashint{%
 \operatorname%
 {\,\,\text{\bf--}\kern-.98em\DOTSI\intop\ilimits@\!\!}}
 \makeatother

\newcommand\bR{\mathbb{R}}

\newcommand\bS{\mathbb{S}}

\newcommand\cR{\mathcal{R}}

\newcommand\cL{\mathcal{L}}

\newcommand\sfu{{\sf u}} 
\newcommand\sfv{{\sf v}}

\newcommand\infsup{\operatornamewithlimits{inf\,\,\,sup}}

\newcommand\loc{\textnormal{loc}}
\newcommand\esssup{\operatornamewithlimits{ess\,sup\,}}
\newcommand{\osc}{\operatornamewithlimits{osc}}
\newcommand\ext{\rm ext}

\begin{document}

\title[$C^{1+\alpha}$-regularity of viscosity solutions]
{$C^{1+\alpha}$-regularity of  viscosity solutions
of general nonlinear parabolic equations}

\author{N.V. Krylov}

\email{nkrylov@umn.edu}
\address{127 Vincent Hall, University of Minnesota,
 Minneapolis, MN, 55455}


\renewcommand{\subjclassname}{%
\textup{2010} Mathematics Subject Classification}


\begin{abstract}
 We investigate the $C^{1+\alpha}$-regularity
of solutions of parabolic equations
$\partial_{t}v+H(v,Dv,D^{2}v,t,x)=0$.
  Our main result  
says that under rather general assumptions
there exist  $C$-viscosity and $L_{p}$-viscosity
solutions which are in $C^{1+\alpha}_{\loc} $.
We allow $H$ to be just measurable in $t$ 
and for its principal part to
have sufficiently small discontinuities
as a function of~$x$. No Lipschitz continuity
of $H$ with respect to $v,Dv$ is required.
\end{abstract}

\maketitle

\mysection{Introduction}
                                        \label{section 3,19,1}

For a real-valued measurable function   $H(\sfu,t,x)$,    \smallskip
$$
\sfu=(\sfu',\sfu''),\quad
\sfu'=\big(\sfu'_{0},\sfu'_{1},...,\sfu'_{d}\big) \in\bR^{d+1},\quad 
\sfu''\in\bS,\quad
(t,x)\in\bR^{d+1}, 
\medskip$$ 
where $\bS$ is the
set of symmetric $d\times d$ matrices, and
sufficiently regular functions $v(t,x)$ we set  \smallskip
$$
H[v](t,x)=H\big(v(t, x),D v(t,x),D^{2}v(t, x),t, x\big),
\medskip$$
and
we will be dealing with 
the parabolic  equations   \smallskip
\begin{equation}
                                                \label{7.29.10}
\partial_t v(t,x)+H[v](t,x)=0
\smallskip\end{equation}

\noindent
in subsets of $[0,T)\times \bR^d $, where $T\in(0,\infty)$
is fixed. Above  \smallskip
$$
\bR^{d}=\big\{x=(x^{1},...,x^{d}):x^{1},...,x^{d}\in \bR\big\},
\medskip$$
$$
\partial_t=\frac{\partial}{\partial t},\quad
 D^{2}u=(D_{ij}u),\quad Du=(D_{i}u),\quad
D_{i}=\frac{\partial}{\partial x^{i}},
\quad D_{ij}=D_{i}D_{j}.
\medskip$$

If $  R\in(0,\infty)$ and $(t,x)\in\bR^{d+1}$, then   \smallskip
$$
B_{R}=\big\{x\in\bR^{d}:|x|<R\big\},\quad B_{R}(x)=x+B_{R},
\medskip$$
$$
C_{ R}=[0,R^{2})\times B_{R},\quad
C_{R}(t,x)=(t,x)+C_{R}.
\medskip$$

We also take a bounded domain $\Omega\subset\bR^{d}$
of class $C^{1,1}$  and set   \smallskip    
$$
\Pi=[0,T)\times\Omega,\quad\partial'\Pi=\bar\Pi 
\setminus\big(\{0\}\times\bar\Omega\big)
\medskip$$

We will be dealing with viscosity solutions of 
\eqref{7.29.10} in $\Pi$. The following definition
is  taken from \cite{CKS00}
and has the same spirit as in  \cite{CIL}.

\begin{definition}
                                             \label{definition 1,29.1}  
  For each
choice of ``regularity'' class $\cR=C $ or $\cR=L_{p}$
we say that $u$ is an $\cR$-viscosity subsolution
\index{$\cR$-viscosity solution}%
 of 
\eqref{7.29.10} in $\Pi$ provided that $u$ is continuous
in $\Pi$ and, for any $\bar C_{r}(t_{0},x_{0})\subset\Pi$  and any function
$\phi$, that is   continuous in $C_{r}(t_{0},x_{0})$ and
whose generalized derivatives satisfy
$\partial_{t}\phi,D\phi,D^{2}\phi\in\cR\big(C_{r}(t_{0},x_{0})\big)$, and is such that $u-\phi$ attains its
maximum over $C_{r}(t_{0},x_{0})$ at $(t_{0},x_{0})$, we have \medskip
\begin{equation}
                                               \label{1,29,2}
 \lim _{ \rho\downarrow0}\esssup_{C_{\rho}(t_{0},x_{0}) }
\big[\partial_{t}\phi(t,x) +
H\big(u(t,x),D\phi(t,x),D^{2}\phi(t,x),t,x\big)\big]\geq 0.
\medskip\end{equation}

In a natural way one defines $\cR$-viscosity supersolution
and calls a function an $\cR$-viscosity solution if it is an
$\cR$-viscosity supersolution and an $\cR$-viscosity subsolution. 
\end{definition}

Note that $C_{r}(t_{0},x_{0})$ contains $\big\{(t,x):
t=t_{0}, |x-x_{0}|<r\big\}$, which is part of its boundary. Therefore,
the conditions like $D^{2}\phi\in C\big(C_{r}(t_{0},x_{0})\big)$ mean that
the second-order derivatives of $\phi$ are continuous
up to this part of the boundary.

\begin{remark}
                              \label{remark 1,29,1}
If $H(\sfu,t,x)$ is a continuous function of $(\sfu,t,x)$
and we are talking about the $C$-viscosity subsolutions,
then \eqref{1,29,2} becomes, of course, just  \medskip
$$
\partial_{t}\phi(t_{0},x_{0}) +
H\big(u(t_{0},x_{0}),D\phi(t_{0},x_{0}),D^{2}\phi( t_{0},x_{0}  ),
 t_{0},x_{0} \big)\geq 0.
\medskip$$
Also note that $L_{p}$-viscosity solutions are
automatically $C$-viscosity solutions.
\end{remark}

The reader is referred 
to \cite{CKS00} for numerous properties of $L_{p}$-viscosity
solutions and to \cite{CIL} for those of $C$-viscosity
solutions.

The notion of $L_{p}$-viscosity solution generalizes that
of $W^{1,2}_{p}$-solution, which is seen from the following
 well-known fact (see \cite{CKS00}).   Set
$$
[\sfu']=(\sfu'_{1},...,\sfu'_{d}).
$$ 
\begin{theorem}
                                          \label{theorem 8,16,1}
Suppose that $H$ is  a  nonincreasing
function of $\sfu'_{0}$, Lipschitz continuous
with respect to $([\sfu'],\sfu'')$ with constant
independent of $\sfu'_{0}$ and $(t,x)$, and 
at all points of its differentiability
with respect to $\sfu''$ we have $D_{\sfu''}H\in\bS_{\delta}$,
where the constant $\delta\in(0,1]$. Suppose that $p\geq d+1$,
$v$ is a continuous in $\bar\Pi$, $L_{p}$-viscosity solution of
\eqref{7.29.10}, and $w\in W^{1,2}_{p,\loc}(\Pi)\cap C(\bar\Pi)$
is a function satisfying \eqref{7.29.10} \(a.e.\) in $\Pi$.
Finally, let $v=w$ on $\partial'\Pi$. Then $v=w$ in $\Pi$.

\end{theorem}

  Our main result  
says that under rather general assumptions
there exist  $C$-viscosity and $L_{p}$-viscosity
solutions which are in $C^{1+\alpha}_{\loc}(\Pi)$.
We allow $H$ to be just measurable in $t$ 
and for its principal part to
have sufficiently small discontinuities
as a function of~$x$.

Wang in Theorem 1.3 of \cite{Wa92_1} assumes a structure 
condition on $H$ which in the case of linear equations
implies that the coefficients of $Du$ are independent 
of $(t,x)$. On the other hand, he proves
the result for {\em any\/} $C$-viscosity solution.

Our Theorem \ref{theorem 12,30.1} contains Proposition 5.4 of \cite{CKS00}
proved there for equations $\partial_{t}u+H(D^{2}u)=0$
apart from the fact that Proposition 5.4 of \cite{CKS00}
is valid for {\em any\/} $C$-viscosity solution.
Theorem \ref{theorem 12,30.1} is also close to Theorem 7.3
of \cite{CKS00}, that is proved for {\em any\/} $L_{p}$-viscosity
solution  and not for the $C$-viscosity or $L_{p}$-viscosity
solution we construct. The most significant
difference in the assumptions is that
$H$ in \cite{CKS00} is supposed to be Lipschitz
in $[\sfu']:=(\sfu'_{1},...,\sfu'_{d})$ and for any $\sfu$ be uniformly in $(t,x)$
close to a function that is uniformly continuous in $(t,x)$
(and not only in $x$).

In \cite{ST_16} the $C^{1+\alpha}$-regularity  
is investigated when $\bar G$ is summable to different powers 
in $t$ and $x$. Again $H$ is Lipschitz in $\sfu$
and satisfies the continuity condition like
  in $(t,x)$ like in \cite{CKS00}.

In what concerns the fully nonlinear {\em elliptic\/} equations,  
Caffarelli \cite{Caf89}
 and Trudinger \cite{Tr_88}, \cite{Tr_89}
were the first authors who proved $C^{1+\alpha}$
regularity for $C$-viscosity solutions of equations  \smallskip
\begin{equation}
                                              \label{4.4.1}
H[u]=f
\end{equation}
without convexity assumptions on $H$. The assumptions in these 
papers are different. 
In \cite{Caf89} the function $H(\sfu,x)$
 is independent of $\sfu'$  
and, for each $\sfu''$, is uniformly sufficiently close to a 
function which is
 continuous
with respect to $x$. In \cite{Tr_88} and 
\cite{Tr_89} the function $H$ depends on all arguments but is
H\"older continuous in $x$. If we ignore the difference
between $C$-viscosity and $L_{p}$-viscosity
solutions, the next step in what concerns
$C^{1+\alpha}$-estimates   was done by   
\'Swi{\c e}ch \cite{Sw_97},
 who considered general $H$ and imposed
the same condition as in \cite{Caf89} on the $x$-dependence,
which is  much weaker than in \cite{Tr_88} and \cite{Tr_89} (but
also imposed the Lipschitz condition on the dependence of $H$ on
$\sfu'$,
whereas
 in \cite{Tr_88} and \cite{Tr_89} only the 
continuity with respect to
$\sfu' $ is assumed). It is worth emphasizing
that these results are about {\em any\/} continuous viscosity
solution and not only about the ones we construct.
 The same bears on the results in \cite{SS_2014},
where the {\em boundary\/} $C^{1+\alpha}$ regularity is obtained.

The best value of $\alpha$ is largely unknown. However, it is proved
in \cite{ST_2015} in case $H=H(\sfu'')$ that the solutions
of $H[u]=0$ are almost $C^{1+\alpha}$ regular if
the solutions of $\hat H[u]=0$ are $C^{1+\alpha}$ regular,
where   
$$
\hat H( \sfu'' ):=\lim_{\tau\to\infty}
\frac{1}{\tau}H ( \tau \sfu''  )
\smallskip$$
assuming that the limit exists.

By the way, it follows from a local version of
Theorem 2.1 of \cite{Kr_13.1}
that, if $\hat H$ is convex, then 
the solutions
of $H[u]=0$ are in
$  W^{2}_{p,\loc}$
for any $p>1$, and thus in $C^{1+\alpha}$
for any $\alpha<1$. This covers Corollary 1.2
of \cite{ST_2015}.

 \mysection{Main results}

                                                \label{section 7,7,1}

Recall some definitions. For $\kappa\in(0,1]$ and 
functions $\phi(t,x)$ on $ \Pi $ set    \smallskip
$$
[\phi]_{C^{ \kappa}( \Pi )}=
\sup_{(t,x),(s,y)\in
 \Pi }\frac{\big|\phi(t,x)-\phi(s,y)\big|}
{|t-s|^{\kappa/2}+|x-y|^{\kappa}},\quad
\|\phi\|_{C ( \Pi )}=\sup_{ \Pi }|\phi|,
\medskip$$  \smallskip
$$
\|\phi\|_{C^{ \kappa}( \Pi )}=\|\phi\|_{C ( \Pi )}
+[\phi]_{C^{ \kappa}( \Pi )}.
\medskip$$
For $\kappa\in(1,2]$ and sufficiently
\index{$N$@Norms!$\vert\vert\phi\vert\vert_{C^{\kappa}(\Pi)}$}%
\index{$A$@Sets of functions!$C^{\kappa}(\Pi)$}%
 regular $\phi$ set
\smallskip
$$
[\phi]_{C^{ \kappa}(  \Pi )}=
\sup_{t,s\in[0,T), x\in\Omega }\frac{\big|\phi(t,x)-\phi(s,x)\big|}
{|t-s|^{\kappa/2} }
\medskip$$ \smallskip
$$
+\sup_{x,y\in \Omega  ,t\in[0,T)}\frac{\big|D\phi(t,x)-D\phi(t,y)\big|}
{ |x-y|^{\kappa-1}},\quad
\|\phi\|_{C^{ \kappa}(  \Pi )}=
\|\phi\|_{C^{1}( \Pi )}+[\phi]_{C^{ \kappa}(  \Pi )}.
\medskip$$

   The set of functions
with finite norm $\|\cdot\|_{C^{\kappa}( \Pi )}$ 
is denoted by $C^{\kappa}( \Pi )$.
Observe that any $u\in C^{\kappa}( \Pi )$ admits a unique
extension to $\bar\Pi$ and we will always consider
this extension while dealing with $u(T,x)$.
For $\kappa=2$ we prefer a less ambigues notation 
$W^{1,2}_{\infty}(\Pi)$ instead of $C^{2}(\Pi)$.
As usual 
\index{$A$@Sets of functions!$C^{\kappa}_{\loc}(\Pi)$}%
we write $u\in C^{\kappa}_{\loc}( \Pi )$
if $u\in C^{\kappa}\big(C_{R}(t,x)\big)$ for any $C_{R}(t,x)$
such that $\bar C_{R}(t,x)\subset \Pi$.

For $\kappa\in (1,2)$ we are also going to use the spaces $C^{1+\kappa}
(\Pi)$ of functions $u\in W^{1,2}_{\infty}(\Pi)$ with finite norm
$$
\|u\|_{C^{1+\kappa}(\Pi)}
=\|v\|_{C^{2}(  \Pi)}+[u]_{C^{1+\kappa}(\Pi)},
$$
where
$$
[u]_{C^{1+\kappa}(\Pi)}=[\partial_{t}u]_{C^{ \kappa-1}(\Pi)}
+[D^{2}u]_{C^{ \kappa-1}(\Pi)}.
\medskip$$

\begin{remark}
                                  \label{remark 12,14.2}
Sometimes it is useful to invoke the well-known embedding theorem
according to which, if $\kappa\in(1,2)$
and $\phi\in C^{\kappa}(C_{R})$, then  \smallskip
\begin{equation}
                                                    \label{12,14.3}
 \big|D\phi(t,x)-D\phi(s,x)\big| 
 \le N(d )|t-s|^{(\kappa-1)/2}[\phi]_{C^{\kappa}(C_{R})}
\smallskip\end{equation}

\noindent whenever $(t,x),(s,x)\in C_{R}$.

For $\delta\in(0,1]$ denote
$$
\bS_{\delta}=\{a\in\bS:\delta^{-1}|\lambda|^{2}
\geq a^{ij}\lambda^{i}\lambda^{j}\geq \delta|\lambda|^{2},
\forall \lambda\in\bR^{d}\}.
$$

\end{remark}

\begin{assumption}
                                    \label{assumption 12,29.1}
 (i) The   function  $H(\sfu,t,x)$
is  
  Lipschitz continuous in $\sfu''$ for every $\sfu',(t,x)\in\bR^{d+1}$
and  at all points of differentiability of 
$H(\sfu,t,x)$  with
respect to $\sfu''$, we have
 $
D_{\sfu''}H  \in \bS_{\delta}
 $, where $\delta$ is a fixed constant in $(0,1]$.

(ii) Either

(a) for any $\sfu_{0}\in\bR^{d+1}\times\bS$
and $(t_{0},x_{0})\in\Pi$ the function $H(\sfu,t,x)$
is continuous at $\sfu_{0}$ as a function of $\sfu$
uniformly with respect to $(t,x)$ belonging
to a neighborhood of $(t_{0},x_{0})$;

\noindent
or

(b)  For any $M\in(0,\infty)$  the function $H(\sfu,t,x)$ is continuous
with respect to $[\sfu']=\big(\sfu'_{1},...,\sfu'_{d}\big)$ uniformly 
with respect to $|\sfu'_{0}|\leq M$, $u''\in\bS$, $(t,x)\in\Pi$.

\end{assumption}

We fix some constants   \smallskip
$$
p\in(d+2,\infty),\quad
K_{0},K_{1},\in[0,\infty),\quad
 R_{0}\in(0,1],\quad\hat R_{0}\in
[0,\infty],
\medskip$$ 
a nonnegative function
$$
\bar G\in L_{p}(\Pi),
\medskip$$
and
 in the following assumption also use
$\theta_{0}=\theta_{0}( d,\delta)\in(0,1]$, whose value  
is  specified in the proof of Lemma \ref{lemma 1,3.4}.

 \begin{assumption}
                                  \label{assumption 12,29.2}
We have a representation  \smallskip
$$
H(\sfu,t,x)=F\big(\sfu'_{0},\sfu'',t,x\big)+G(\sfu,t,x).
$$ 

(i) The functions $F$ and $G$ are measurable functions 
of their arguments,  {\em continuous\/}  with respect to $\sfu$
for any $(t,x)$.

(ii) For all values of the arguments   \smallskip
$$
\big|G(\sfu,t,x)\big|\leq K_{0}|\sfu'|+\bar G(t,x) .
$$
 
(iii)
 The function $F$ is Lipschitz continuous with respect to $\sfu''$
and at all points   of differentiability of $F$
with respect to $\sfu''$
 and any
$\varepsilon\in\bS$ such that $|\varepsilon|=1$  we have  \smallskip
$$
D_{\sfu''}F
 + \theta_{0} \varepsilon \in\bS_{\delta/2}.
\medskip$$

(iv) For any $\sfu'_{0},\sfv'_{0} \in\bR$, $(t ,x_{i})\in\bR^{d+1}$,
$i=1,2$, and $\sfu''\in\bS$ we have   \smallskip
$$
\big|F\big(\sfu'_{0},\sfu'',t ,x_{1}\big)-F\big(\sfv'_{0},\sfu'',t 
,x_{2}\big)\big|\leq
K_{1}+\theta_{0}|\sfu''| 
\medskip$$
as long as $|x_{1}-x_{2}|\leq R_{0}$
and $\big|\sfu'_{0}-\sfv'_{0}\big|<\hat  R_{0}$.

(v) We have $F\big(\sfu'_{0},0,t,x\big)=0$ for all
$\sfu'_{0} \in\bR$, $(t ,x )\in\bR^{d+1}$.

\end{assumption}

\begin{remark}
                                               \label{remark 1,25,4}
The value $\infty$ is allowed for $\hat R_{0}$
with the purpose to cover the cases in which $F$
is almost independent of $\sfu'_{0}$ when conditions
like \eqref{1.25.7} below are automatically satisfied.
\end{remark}

\begin{remark}
                                               \label{remark 3,22,2}
 Observe that Assumption \ref{assumption 12,29.2} (iv)
is automatically satisfied if $F$ is a function of only
$\sfu''$ and $t$.
\end{remark}

\begin{assumption}
                      \label{assumption 12,30.3}
  We are given
$g\in C(\bar\Pi) $.
 
\end{assumption}  

\begin{theorem}
                                    \label{theorem 12,30.1}
There is a constant $\theta_{0}=\theta_{0}(d,\delta)\in(0,1]$
such that, if the above assumptions are satisfied
with this $\theta_{0}$, then  there exist a 
$\kappa=\kappa(d,\delta,p)\in(1,2)$ and a function $v\in
 C^{\kappa}_{\loc}(\Pi)\cap C(\bar\Pi) $
that  is a  $C$-viscosity or an $L_{p}$-viscosity solution of 
  equation  \eqref{7.29.10}
in $\Pi$   with boundary 
condition $v =g $ on $\partial'\Pi$
according as requirements \(a\) or \(b\)
in Assumption \ref{assumption 12,29.1} \(ii\,\) are satisfied.
 
Furthermore, for any $r,R\in(0,R_{0}]$ satisfying $r<R$ and $(t,x)
\in\Pi$ such that $C_{R}(t,x)\subset \Pi$ and
\begin{equation}
                                                 \label{1.25.7}
\osc_{C_{R}(t,x)}v<\hat R_{0}
\end{equation}
 we have
\begin{equation}
                                                 \label{1.25.6}
[v]_{C^{\kappa}(C_{r}(t,x))}
\leq N(R-r)^{-\kappa}\sup_{C_{R}(t,x)}|v|
+N \big(K_{1}+\|\bar G\|_{L_{p}(C_{R}(t,x))}\big),
\end{equation}
where $N$ depend  only on $d,\delta$, and $K_{0}$.

\end{theorem}

This theorem is proved in Section \ref{section 1,3.1}.

 \begin{remark}
We assumed that $\Omega\in C^{1,1}$ just to be able to refer
to the results available at this moment, but actually
much less is needed for Theorem \ref{theorem 12,30.1}
to hold. For instance the exterior cone condition would
suffice.

\end{remark}

\begin{example}
                                          \label{example 1,31,1}
Let $A$ and $B$ be countable sets and suppose that
for any $\alpha\in A$ and $\beta\in B$ we are given
an $\bS_{\delta}$-valued functions 
$a^{\alpha\beta}\big(\sfu'_{0},t,x\big)$ and 
a
real-valued
function $b^{\alpha\beta}(\sfu',t,x)$
defined for $\sfu',(t,x)\in\bR^{d+1}$. Assume that
these functions are measurable as functions of $(\sfu',t,x)$,
continuous with respect to $ \sfu' $
  uniformly with respect
to $\big(\alpha,\beta,(t,x)\big)\in A\times B\times\Pi$, and 
$a^{\alpha\beta}\big(\sfu'_{0},t,x\big)$ is continuous with respect to
$x$ uniformly with respect to 
$\big(\alpha,\beta,\sfu'_{0},t\big)\in A\times B\times\bR\times\bR$. Finally,
suppose that
for all values of indices and arguments
$$
\big|b^{\alpha\beta}(\sfu',t,x)\big|\leq K_{0}|\sfu'|+\bar G(t,x) .
$$
Then the following equation
$$
\partial_{t}v+\infsup_{\alpha\in A\,\,\beta\in B}
\Big[a^{\alpha\beta}_{ij}(v,t,x)D_{ij}v+b^{\alpha\beta}
(v,Dv,t,x)\Big]=0
$$
in $\Pi$ with boundary condition $v=g$ on $\partial'\Pi$
has an $L_{p}$-viscosity solution which belongs to
$ C^{\kappa}_{\loc}(\Pi)\cap C(\bar\Pi) $.

This follows immediately from Theorem \ref{theorem 12,30.1}
if one sets
$$
F\big(\sfu'_{0},\sfu'',t,x\big)=
\infsup_{\alpha\in A\,\,\beta\in B}
 a^{\alpha\beta}_{ij}\big(\sfu'_{0},t,x\big)\sfu''_{ij}
$$
and observes that, for any $\theta_{0}$, in particular,
for the one from Theorem \ref{theorem 12,30.1},
one can find $R_{0}$ and $\hat R_{0}$, for which
Assumption \ref{assumption 12,29.2} (iv)
is satisfied with an appropriate $K_{1}$.

One can also see that the continuity of 
$a^{\alpha\beta}\big(\sfu'_{0},t,x\big)$ with respect to
$x$ can be relaxed allowing sufficiently small
discontinuities. It is also worth noting that
in \cite{CKS00} in case of the Isaacs equations
$a$ is independent of $\sfu'_{0}$ and $b$ is an affine function
of $\sfu'$.

\end{example}

We explain the main rough ideas in the proof of Theorem \ref{theorem 12,30.1}
in case $F$ is independent of $\sfu'_{0}$ and $\bar G$ is bounded.
It
consists of establishing  a priori estimates of the type
\begin{equation}
                                                  \label{7,7,1}
\sup_{C_{r}(t_{0},x_{0})}|v-l|\leq Nr^{\kappa}
\end{equation}
for any $C_{r}(t_{0},x_{0})$ which is strictly inside $\Pi$
and any small $r>0$ as long as an affine function $l=l(x)$
is chosen appropriately. This turns out to be enough
to get an estimate of the $ C^{\kappa}$-norm of $v$ in small cylinders
which are strictly inside $\Pi$ (see Lemma \ref{lemma 10.23.1}).
Then,
to obtain \eqref{7,7,1} we represent $v$ as $h+w$, where
$h$ is a solution of 
\begin{equation}
                                                  \label{7,7,2}
\partial_{t}h+F[h]=0
\end{equation}
in $C_{r}(t_{0},x_{0})$ with boundary data $v$ and
$w=v-h$ is found from
$$
0=\partial_{t}v+H[v]-\partial_{t}h-F[h]
=\partial_{t}w+a^{ij}D_{ij}w+G[v],
$$
where $(a^{ij})$ is a certain $\bS_{\delta}$-valued function.
Since $w=0$ on $\partial'C_{r}(t_{0},x_{0})$ by the maximum
principle $|w|\leq Nr^{2}\big(1+\sup(|Dv|,C_{r}(t_{0},x_{0})\big)$.
The latter supremum is irrelevant because we are going
to estimate the $ C^{\kappa}$-norm of $v$ and $\kappa>1$.

 Then we see that, to get \eqref{7,7,1}, it suffices to
prove that \eqref{7,7,1} holds with $h$ in place of $v$.
To do this step we want to replace $F$ with the one
independent of $(t,x)$. Freezing the coefficients
does not help because there is no hope to control
the second-order derivatives of solutions of such equations.
Therefore, we just replace $F$ with 
$$
F_{0}^{(\pm)}=F(\sfu'',t,x_{0})\pm\big(K_{1}+\theta_{0}|\sfu''|\big)
$$
and introduce $v^{(\pm)}$ as solutions of
\begin{equation}
                                                  \label{7,7,3}
\partial_{t}v^{(+)}+F_{0}^{(+)}[v^{(+)}]=0,
\end{equation}
\begin{equation}
                                                  \label{7,7,4}
\partial_{t}v^{(-)}+F_{0}^{(-)}[v^{(-)}]=0
\end{equation}
in $C_{r}(t_{0},x_{0})$ with boundary condition $v$.
Since $F_{0}^{(-)}\leq F\leq F_{0}^{(+)}$ in $C_{r}(t_{0},x_{0})$
if $r$ is small enough, we have 
  $v^{(-)}\leq h\leq v^{(+)}$ by the maximum principle.
Furthermore, since $F_{0}^{(\pm)}$ is independent of $x$,
one can differentiate the equations
\eqref{7,7,3} and \eqref{7,7,4} with respect to $x$ and get
  estimates of the H\"older constants of $Dv^{(\pm)}$
by the Krylov-Safonov theorem.
These estimates guarantee that $v^{(\pm)}$ can be approximated
by affine functions of $x$ as in \eqref{7,7,1}.

Then the only thing which remains is to estimate $|v^{(+)}
-v^{(-)}|$. To this end, for $\theta\in[-\theta_{0},\theta_{0}]$ we introduce 
$$
F_{0} ( \sfu'',t,\theta)=F_{0} ( \sfu'',t)+\theta|\sfu''|
+K_{1}\theta/\theta_{0}
$$
  and
define $v^{\theta}$ from the equation
\begin{equation}
                                                  \label{7,7,5}
\partial_{t}v^{\theta}+F_{0} (D^{2}_{x} v^{\theta},t,\theta)=0
\end{equation}
in $C_{r}(t_{0},x_{0})$ with boundary condition $v$.
Since $v^{\pm\theta_{0}}=v^{(\pm)}$ to estimate
$|v^{(+)}
-v^{(-)}|$, it suffices to estimate $D_{\theta}v^{\theta}$.

Now comes an idea originated in the theory of diffusion
processes. We look at \eqref{7,7,5} as an equation
in variables $(t,x,\theta)$. There is no derivatives
with respect to $\theta$, so that it is a degenerate equation,
but this was never a problem in such matters in that theory,
that suggests that the function
$$
V(t,\theta,x,\tau,\xi):=\tau D_{\theta}v^{\theta}(t ,x)+
\xi^{i}D_{x^{i}}v^{\theta}(t ,x)
$$
satisfies a parabolic equation with respect to
the variables $(t,\theta,x,\tau,\xi)$. Indeed
$$
\partial_{t}V+a^{ij}D_{x^{i}x^{j}}V+\tau\varepsilon^{ij}
D_{x^{i}\xi^{j}}V
+N\tau^{2}\delta^{ij}D_{\xi^{i}\xi^{j}}V+K_{1}\tau/\theta_{0}=0,
$$
where $(a^{ij})$ is an $\bS_{\delta/2}$-valued function
(see Assumption \ref{assumption 12,29.2} (iii)),
$$
\varepsilon^{ij}=D_{x^{i}x^{j}}v^{\theta}/|D^{2}_{x}v^{\theta}|,
$$ 
and $N$ is any constant. Since there is no derivatives with respect to $\tau$,
it is just a parameter and for $W(t,\theta,x,\xi):=V(t,\theta,x,1,\xi)$
we obtain the equation
$$
\partial_{t}W+a^{ij}D_{x^{i}x^{j}}W+ \varepsilon^{ij}
D_{x^{i}\xi^{j}}W
+N \delta^{ij}D_{\xi^{i}\xi^{j}}W+K_{1} /\theta_{0}=0,
$$
which is parabolic if $N$ is sufficiently large. We consider
this equation in $(t,x,\xi)\in C_{r}(t_{0},x_{0})\times \bR^{d}$
and see that to estimate $W$, it suffices to have a good control
of it 
  on the parabolic boundary of this set, where 
$D_{\theta}v^{\theta}(t ,x)=0$ by construction. Thus we see
that we need to estimate $|D_{x}v^{\theta}|$ in $ C_{r}(t_{0},x_{0})$.
This will be done by differentiating \eqref{7,7,5}
with respect to $x$ and using the maximum principle,
which reduces the matter to estimating 
$|D_{x}v^{\theta}|$ in $\partial' C_{r}(t_{0},x_{0})$.
 
  As a general comment
we point out that, since we have to have sufficiently
smooth solutions in the above argument, we use cut-off equations
and use finite-differences to avoid using third-order
derivatives.

\mysection{Auxiliary results about 
linear  equations}
                                           \label{section 12,21.1}
We will be using 
a common way of approximating functions
in $ C^{\kappa}(  C_{1})$ by infinitely
differentiable ones.

\begin{lemma}
                                    \label{lemma 12,17.2}
Let $\kappa\in(0, 2) $, $r\in(0,\infty)$, 
$g, h \in C^{\kappa}( C_{ r })$.
Then for any $\varepsilon>0$ there exists an
infinitely differentiable functions 
$g^{\varepsilon}$ and $h^{\varepsilon}$  on $\bR^{d+1}$ such that in $ C_{r}$  \smallskip
$$
|g-g^{\varepsilon}|\leq N[g]_{C^{\kappa}(  C_{r})}
(r\varepsilon)^{\kappa},\quad
|Dg-Dg^{\varepsilon}|\leq N[g]_{C^{\kappa}(  C_{r})}
(r\varepsilon)^{\kappa-1},
\vspace{8pt}$$
\begin{equation}
                                               \label{12,17.5}
 \big|\partial_{t}g^{\varepsilon}\big|+ |D^{2}g^{\varepsilon}|    
+r \varepsilon|D^{3}g^{\varepsilon}|+
r \varepsilon\big|D\partial_{t} g^{\varepsilon}\big|\leq
N[g]_{C^{\kappa}(  C_{r})}(r \varepsilon)^{\kappa-2},
\vspace{10pt}\end{equation} \smallskip
$$
 [  h^{\varepsilon}]_{C^{1+\kappa}( C_{r})}
\leq N[h]_{C^{\kappa}(  C_{r})}(r \varepsilon)^{ -1}.
\medskip$$
\end{lemma}

Proof. Parabolic scalings reduce the general situation to
the one in which $r=1$. Then this well-known result 
 is obtained by first
continuing $g(t,x)$, $h(t,x)$ as   functions of $t$ to $\bR$ 
to become   even,
2-periodic functions, then continuing thus obtained functions
across $|x|=1$ almost preserving $[g]_{C^{\kappa}(  C_{1})}$,
$[h]_{C^{\kappa}(  C_{1})}$
in the whole space
and then taking  convolutions
with $\delta$-like kernels.
The lemma is proved. \qed

\begin{theorem}
                                           \label{theorem 12,19.1}  
 
Let $\kappa\in(1,2)$ and let $g\in C^{\kappa}(  C_{1})$.
Then there exists a unique $u\in C^{\infty}_{\loc}(C_{1})
\cap C^{\kappa}(  C_{1})$ which satisfies the heat equation
\smallskip\begin{equation}
                                                \label{12,19.3}
\partial_{t}u+\Delta u=0
\smallskip\end{equation}
in $C_{1}$ and equals $g$ on $\partial'C_{1}$. Furthermore,
there exists a constant $N=N(d,\kappa)$ such that
\begin{equation}
                                                \label{12,19.1}
[u]_{C^{\kappa}(  C_{1})}\leq N[g]_{C^{\kappa}(  C_{1})},
\medskip\end{equation}
\begin{equation}
                                                \label{12,19.2}
\big|D^{2}u(t,x)\big|\leq N\Big[\big(1-|x|\big)\wedge\sqrt{1-t}\Big]^{\kappa-2}
|g|_{C^{\kappa}(  C_{1})}.
\smallskip\end{equation}
\end{theorem}

Proof. One can subtract an affine function of $x$ from
$g$ and reduce the general situation to the one where
\begin{equation}
                                                \label{12,19.7}
g(1,0)=0,\quad Dg(1,0)=0.
\smallskip\end{equation}
 In that case take $g^{\varepsilon}$ from Lemma \ref{lemma 12,17.2}.
Then by a classical result (see, for instance,
 Theorem 5.14 in \cite{Li} or
Theorems 10.3.3 and 10.2.2 in \cite{Kr_96}) there exists
a unique $u^{\varepsilon}\in   C^{1+\kappa} ( C_{1})$
satisfying \eqref{12,19.3} in $C_{1}$ and equal to $g$
on $\partial'C_{1}$. In addition,
\begin{equation}
                                                \label{12,19.6}
\|u^{\varepsilon}\|_{ C^{1+\kappa} (  C_{1})}
\leq N\|g^{\varepsilon}\|_{ C^{1+\kappa} (  C_{1})}.
\smallskip\end{equation}

Furthermore, by classical results (see, for instance,
 Theorem 8.12.1 in \cite{Kr_96}) the functions
$u^{\varepsilon}$ are infinitely differentiable with respect to $x$ in
$\bar C_{r}$ for any $r<1$ and by Theorem 8.4.4 in \cite{Kr_96}
and the maximum principle
any derivative of $u^{\varepsilon}$ of any order with respect
to $x$  is bounded in $\bar C_{r}$ for any $r<1$ by a constant
independent of $\varepsilon$. As it follows from equation \eqref{12,19.3}
itself, the same holds for derivatives with respect
to $t$ of any derivative of any order with respect to $x$
(cf.~Exercise 8.12.4 in \cite{Kr_96}).

 In addition, by the maximum
principle and \eqref{12,17.5},
\medskip\begin{equation}
                                                \label{12,19.4}
|u^{\varepsilon_{1}}-u^{\varepsilon_{2}}|\leq
|g^{\varepsilon_{1}}-g^{\varepsilon_{2}}|\leq N[g]_{C^{\kappa}( C_{1})}
(\varepsilon_{1}^{\kappa}+\varepsilon_{2}^{\kappa}).
\smallskip\end{equation}

It follows that, as $\varepsilon\downarrow0$, $u^{\varepsilon}$
converges uniformly on $\bar C_{1}$ to a continuous function,
which is equal to $g$ on $\partial'C_{1}$, is infinitely differentiable
in $C_{1}$ and satisfies \eqref{12,19.3}. In light of \eqref{12,19.4}
we have
\smallskip\begin{equation}
                                                \label{12,19.5}
|u -u^{\varepsilon }| \leq N[g]_{C^{\kappa}(  C_{1})}
 \varepsilon ^{\kappa} .
\end{equation}

Also \eqref{12,19.6} and \eqref{12,19.7} along with \eqref{12,17.5}
imply that
\smallskip$$
[u^{\varepsilon}]_{ C^{1+\kappa} (  C_{1})}
\leq N[g^{\varepsilon}]_{ C^{1+\kappa} (  C_{1})}\leq
N[g]_{C^{\kappa}(  C_{1})}
 \varepsilon ^{-1}.
\medskip$$

Now take $x_{0}\in\bR^{d}$, unit $l\in\bR^{d}$ and $0<h_{1}\leq h$
such that
\smallskip$$
x_{0},x_{0}+ h_{1}l_{1},  
x_{0}+2hl,x_{0}+2hl+h_{1}l_{1}\in \bar B_{1}.
\medskip$$
Observe that for any $t\in[0,1]$
\medskip$$
u(t,x)=\big[u(t,x)-u^{\varepsilon}(t,x)\big]+u^{\varepsilon}(t,x_{0})
+(x^{i}-x^{i}_{0})D_{i}u^{\varepsilon}(t,x_{0})
 $$
\medskip$$
+(1/2)(x^{i}-x^{i}_{0})
(x^{j}-x^{j}_{0})D_{ij}u^{\varepsilon}(t,x_{0})+v(t,x),
\medskip$$
where 
\smallskip$$
|v(t,x)|\leq N|x-x_{0}|^{1+\kappa}
[D^{2}u^{\varepsilon}]_{C^{\kappa-1}(C_{1})}
\leq N|x-x_{0}|^{1+\kappa}[g]_{C^{\kappa}(  C_{1})}
 \varepsilon ^{-1}.
\medskip$$
Since the third-order finite difference of any quadratic
polynomial is zero and
$$
\big|u(t,x)-u^{\varepsilon}(t,x)\big|\leq N[g]_{C^{\kappa}(  C_{1})}
 \varepsilon ^{\kappa},
$$ 
 we have   \smallskip
$$
\big|(T_{h ,l }-1)^{2}(T_{h_{1},l_{1}}-1) 
u(t,x_{0})\big|\leq N[g]_{C^{\kappa}( C_{1})}
 \varepsilon ^{\kappa} 
+Nh^{1+\kappa}[g]_{C^{\kappa}(  C_{1})}
 \varepsilon ^{-1}.
\medskip$$
By taking $\varepsilon=h$
  we arrive at  \smallskip
$$
\big|(T_{h ,l }-1)^{2}(T_{h_{1},l_{1}}-1) 
u(t,x_{0})\big|\leq N[g]_{C^{\kappa}(  C_{1})} h^{\kappa}.
\medskip$$

As it can be shown (or extracted from \cite{Go_62}) that, the arbitrariness
of $x_{0}$, $h$, $h_{1}$, $l$, $l_{1}$ in the above inequality implies that
 for any $t\in[0,1]$   \smallskip
$$
\big[u(t,\cdot)\big]_{C^{\kappa}(  B_{1})}
\leq N\big([g]_{C^{\kappa}(  C_{1})}+\osc_{C_{1}}u\big),
$$
which along with  the maximum principle    
show that
$$
\big[u(t,\cdot)\big]_{C^{\kappa}(  B_{1})}
\leq N\big([g]_{C^{\kappa}(  C_{1})}+\osc_{\partial' C_{1}}g\big),
\medskip$$
where $\osc_{\partial' C_{1}}g$ can be replaced with 
$[g]_{C^{\kappa}( C_{1})}$ in light of \eqref{12,19.7}.

Next, fix $x\in B_{1}$ and take $t_{0}\in(0,1)$, $h>0$,
such that $t_{0}+2h^{2}\in(0,1)$. Observe that   \smallskip
$$
u(t,x)=\big[u(t,x)-u^{\varepsilon}(t,x)\big]+u^{\varepsilon}(t_{0},x)
+(t-t_{0})\partial_{t}u^{\varepsilon}(t_{0},x)+w(t,x),
\medskip$$
where
$$
\big|w(t,x)\big|\leq|t-t_{0}|^{(\kappa+1)/2}
[\partial_{t}u^{\varepsilon}]_{C^{\kappa-1}(C_{1})}
\leq N|t-t_{0}|^{(\kappa+1)/2}
[g]_{C^{\kappa }(C_{1})}\varepsilon^{-1}.
\medskip$$
Since the second-order differences of linear function
are equal to zero,  \smallskip
$$
\big|u(t_{0}+2h^{2},x)-2u(t_{0}+h^{2},x)+u(t_{0},x)\big|
\medskip$$
$$
\leq  N[g]_{C^{\kappa}(  C_{1})}
 \varepsilon ^{\kappa}+Nh^{1+\kappa}[g]_{C^{\kappa}(  C_{1})}
 \varepsilon ^{-1}.
\medskip$$ 
Here, for $\varepsilon=h$, the right-hand side becomes  \smallskip
$$
N[g]_{C^{\kappa}( C_{1})}
 h^{\kappa}
\medskip$$
which implies that  \smallskip
$$
\big|u(t,x)-u(s,x)\big|\leq  |t-s|^{\kappa/2}
 N\big([g]_{C^{\kappa}(  C_{1})}+\osc_{(0,1)}u(\cdot,x)\big) 
\medskip$$
if $s,t\in[0,1]$.
Again the last oscillation can be replaced by $[g]_{C^{\kappa}(  C_{1})}$,
and this proves \eqref{12,19.1}.

By Theorem 8.4.4 in \cite{Kr_96} for any $(t_{0},x_{0})\in C_{1}$  \smallskip
$$
\big|\partial_{t}u(t_{0},x_{0})\big|+\big|D^{2}u(t_{0},x_{0})\big|
\leq NR^{-2}\sup_{C_{R}(t_{0},x_{0})}|u|,
\medskip$$
where $R=  (1-|x_{0}|)\wedge\sqrt{1-t_{0}} $.
This holds for any sufficiently
regular solution of \eqref{12,19.3} and not only for the one constructed above.
 Therefore $u$
on the right can be replaced with $u-\hat u$, where
$\hat u$
is any affine function of $x$. By taking $\hat u$
as the first-order Taylor polynomial of $u(t_{0},x)$ with respect to
$x$ at $ x_{0} $ and using \eqref{12,19.1} we come to
  \eqref{12,19.2}. The theorem is proved. \qed

\begin{lemma}
                                           \label{lemma 12,26.10} 
For $\kappa\in(1,2)$ there is a function $\Phi\in C^{1,2}_{\loc}(C_{1})
\cap C(\bar C_{1})$ and a constant $N=N(\kappa,\delta)$ such that 
\begin{equation}
                                              \label{12,27.1}
\partial_{t}\Phi+a^{ij}D_{ij}\Phi\leq-
\Big[\big(1-|x|\big)\wedge\sqrt{1-t}\Big]^{\kappa-2}
\end{equation}
in $C_{1}$ for any $(a^{ij})\in\bS_{\delta}$ and   
$$
0\leq \Phi(t,x)\leq N \big(1-|x|\big) .
 $$
\end{lemma}

Proof. Set $\rho(x)= 1-|x|^{2} $, $\beta=(\kappa+1)/2$,
 and, for a constant $N_{0}$
to be determined later, define   \smallskip
$$
\Phi (t,x)=N_{0}\big[\rho(x)-(1/\kappa)\rho^{\kappa}(x)\big]
+\beta^{-1}(1-t)^{\beta}\rho(x).
\smallskip$$
We have
$$
D_{i}\rho=-2x^{i},\quad D_{ij}\rho=-2\delta^{ij},
\medskip$$
$$
D_{ij}\Phi=-2N_{0}[1-\rho^{\kappa-1}]\delta^{ij}
-4N_{0}(\kappa-1)\rho^{\kappa-2}x^{i}x^{j}-2\beta^{-1}
(1-t)^{\beta}\delta^{ij}.
\medskip$$
Hence,
$$
\partial_{t}\Phi+a^{ij}D_{ij}\Phi
\leq-\delta I\rho^{\kappa-2} 
- (1-t)^{\beta-1}\rho,
\medskip$$
where
$$
I:=2N_{0}[\rho^{2-\kappa}-\rho]+4N_{0}(\kappa-1)|x|^{2}.
\medskip$$

Obviously, there exists $N_{0}=N_{0}(\delta,\kappa)$ such that in $C_{1}$
we have $
I\geq 2/\delta 
 $,
in which case  \smallskip
$$
-\delta I\rho^{\kappa-2} \leq-2\rho^{\kappa-2}
=-2(1+|x|)^{\kappa-2}(1-|x|)^{\kappa-2}\leq-(1-|x|)^{\kappa-2}
\medskip$$
and \eqref{12,27.1} holds if $1-|x|\leq \sqrt{1-t}$. In case
$1-|x|\geq \sqrt{1-t}$ we have $\rho\geq \sqrt{1-t}$ and
$$
- (1-t)^{\beta-1}\rho\leq- (1-t)^{\beta-1/2}=-(1-t)^{\kappa/2-1},
\medskip$$
so that \eqref{12,27.1} holds again. The lemma is proved.    \qed

\mysection
{Estimating the difference of solutions
of two different equations}
                                                  \label{section 12,31.1}

We need the following Theorem 6.1 of \cite{Kr_17.1}.

\begin{theorem}
                                        \label{theorem 4,5,1}  

  Suppose that Assumption \ref{assumption 12,29.1} \(i\,\)
is satisfied, $\Omega\in C^{1,1}$, $H$ is a continuous function 
of $\sfu$,  the number
$$
\bar{H} :=\sup_{\sfu', t, x}\big(|H (\sfu',0,t, x)|-K_{0}|\sfu'|\big)
\quad(\geq0)
$$
is finite, and
$g\in W^{1,2}_{\infty }(\bR^{d+1})$.
Then there exists a convex positive homogeneous of degree
one function $P(\sfu'')$ such that 
   at all points of its differentiability 
$D_{\sfu''}P \in
\bS_{\bar\delta}$, where $\bar\delta=\bar\delta(d,\delta)\in(0,\delta)$,
and for $P[u]=P(D^{2}u)$
and any $K>0$ the equation
\begin{equation}
                                               \label{9.28.4}
\partial_{t}v+ \max(H[v],P[v]-K)=0  
\end{equation}
in $\Pi$ with boundary condition $v=g$ on $\partial'\Pi$
has a solution  $v\in
 W^{1,2}_{p}(\Pi)$ for any $p\geq 1$.
\end{theorem}

Let $F_{0}(\sfu'',t)$ be a function satisfying   
Assumption \ref{assumption 12,29.2}  (iii), measurable in $t$,
and such that $F_{0}(0,t)=0$. 
Fix a constant $K_{1}\geq0$ and
set    \smallskip
$$
F^{(\pm)}(\sfu'',t)=F_{0}(\sfu'',t )\pm\big(K_{1}+\theta_{0}|\sfu''|\big),
\medskip$$
take $P(\sfu'')$ from Theorem \ref{theorem 4,5,1}
 with $\delta/2$ in place of $\delta$ 
and for   fixed $R\in(0, \infty) ,K>0$   
 consider the equations   \smallskip
$$
\partial_{t}v+ \max\big(F^{(\pm)}[v], P[v]-K\big)=0
\medskip$$
in $C_{R}$ with boundary data $v=g$ on $\partial'C_{R}$,
where $g\in W^{1,2}_{\infty}(
C_{R})$ is a given function. By $v^{(\pm)}$
we denote their solutions that exist by Theorem \ref{theorem 4,5,1}
and belong to $W^{1,2}_{p}(C_{R})$ for any $p\in[1,\infty)$.

\begin{theorem}
                                          \label{theorem 12,28.1}

For any $\kappa\in(1,2)$ there exists a constant
$N=N(\kappa,\delta,d)$ such that in $\bar C_{R}$ we have  \smallskip
$$  
\big|v^{(+)}-v^{(-)}\big|\leq NR^{2} K_{1} +NR^{\kappa}\theta_{0}
[g]_{C^{\kappa}(C_{R})}.
\smallskip$$
\end{theorem}

The proof of this theorem is based on the following
two auxiliary results.
\begin{lemma}
                                       \label{lemma 1.1.5}
Let $a$ be a measurable $\bS_{\delta}$-valued function, $p>d+2$,
and let $v\in W^{1,2}_{p}(C_{1})$ be a solution
of   \smallskip
$$
\partial_{t}v+a^{ij}D_{ij}v+f=0
\medskip$$
in $C_{1}$ \(a.e.\)  with boundary condition $v=g$
on $\partial'C_{1}$, where $|f|\leq\bar f$
for a constant $\bar f\in[0,\infty)$.
Assume that \eqref{12,19.7} holds. Then,
for any $\kappa\in(1,2)$,
there exists $N=N(\delta,d,\kappa)$ such that
\begin{equation}
                                                \label{12,27.30}
\sup_{(0,1)\times\partial B_{1}}|D v |
\leq N [Dg]_{C^{\kappa}(C_{1})}+N \bar f.
\end{equation}
\end{lemma}

Proof. By the embedding 
Lemma 2.3.3 of \cite{LSU}, the function $D v$
is continuous in $\bar C_{1}$.
Then take the function $u$ from
 Theorem \ref{theorem 12,19.1}
and set $w(t,x)=v( t,x)-u(t,x)$. We have  \smallskip
$$
 \partial_{t}w+a^{ij}D_{ij}w+h=0,
\medskip$$
where   $h:=f+ \partial_{t}u+a^{ij}D_{ij}u$.
By using  
Lemma \ref{lemma 12,26.10}
 we conclude    \smallskip
$$
\big|\partial_{t}w+a^{ij}D_{ij}w\big|\leq  \bar f +
N \Big[\big(1-|x|\big)\wedge\sqrt{1-t}\Big]^{\kappa-2}
|g|_{C^{\kappa}(  C_{1})}
\medskip$$
$$
\leq - 
 \partial_{t}\Psi-a^{ij}D_{ij}\Psi,
\medskip$$
where $\Psi=N|g|_{C^{\kappa}(  C_{1})}\Phi+ \bar f \delta^{-1}
\big(1-|x|^{2}\big)$.
By the maximum principle    \smallskip
$$
\big|w(t,x)\big|\leq \Psi (t,x)
\leq N_{0}\big[|g|_{C^{\kappa}(C_{1})}+ \bar f \,\big]\big(1-|x|\big).
\medskip$$
Consequently, the normal derivative of $v( t,x)$
at a point $x_{0}\in \partial B_{1}$  by magnitude
is less than the absolute value of the normal derivative
of $u(t,x)$ plus $N_{0}\big(|g|_{C^{\kappa}(C_{1})}+ \bar f\, \big)$. By
interpolation inequalities $|Du|$ is estimated by
$[u]_{C^{\kappa}(C_{1})}$ and $\osc_{C_{1}}u\leq\osc_{C_{1}}g$,
where the former is estimated in \eqref{12,19.1} by 
$|g|_{C^{\kappa}(C_{1})}$ and the latter is estimated by the same
quantity due to \eqref{12,19.7}. Thus, the normal derivative
of $v( t,x)$ admits the estimate we are after.
By noting that the tangential derivatives of $v( t,x)$
coincide with those of $g(t,x)$, we finally come to
\eqref{12,27.30}. The lemma is proved.   \qed

\begin{lemma}
                                    \label{lemma 12,18.1}
Let $ R, \chi\in(0,\infty)$ and $\kappa\in(1,2)$ be  constants 
and let $F(\theta,\sfu'',t)$ be a measurable
function on $[-\theta_{0},\theta_{0}]\times
\bS\times\bR$, which is   Lipschitz continuous
in $(\theta,\sfu'')$ for any $t$ and such that 
 at all points of its differentiability \smallskip
$$
|D_{\theta}F |\leq \chi+|\sfu''|,\quad D_{\sfu''}F\in\bS_{\delta}.
\smallskip$$
Also suppose that $F(\theta,0,t)$ is bounded.
Take   $g\in W^{1,2}_{\infty}(
C_{R}) $, and assume that
for any $\theta\in [-\theta_{0},\theta_{0}]$ the equation  \smallskip
\begin{equation}
                                                \label{12,18.1}
\partial_{t}v+F(\theta,D^{2}v,t)=0
\medskip\end{equation}
in $C_{R}$ \(a.e.\) with boundary condition $v=g$ on $\partial' C_{R}$
has a solution $v=v(\theta,\cdot)\in W^{1,2}_{p}(
C_{R}) $, where $p>d+2$, $p\geq 2d+1$. Then  for any $(t,x)\in C_{R}$
the function $v(\theta,t,x)$ is Lipschitz continuous
on $[-\theta_{0},\theta_{0}]$ and at all points of its differentiability
with respect to $\theta$  \smallskip
\begin{equation}
                                                \label{12,18.2}
|D_{\theta}v|\leq  R^{2}\big(\chi   +N\sup_{t}|F(\theta,0,t)|\big) 
 +NR^{\kappa}
[g]_{C^{\kappa}(C_{R})},
\medskip\end{equation}
where the constant $N$ depends only on $\delta,d,\kappa$.
\end{lemma}

Proof. The idea of the proof comes from the theory
of diffusion processes 
and is explained at the end of Section \ref{section 7,7,1}.
The parabolic equation for the directional derivative of $v(\theta,t,x)$
with respect to $(\theta,x)$ in appropriate directions
will be what we are interested in.
Since in our situation there is no guarantee that
$v$ is smooth enough, we follow Trudinger's method (see \cite{Tr_88})
 based on finite-differences.
 
As usual, parabolic scalings reduce the general situation
to the one in which $R=1$. Also by subtracting   from 
$g$ and $v$ the same affine  function of $x$ we may assume that
\eqref{12,19.7} holds.
  
 In that case
fix $\theta^{0}\in(-\theta_{0},\theta_{0})$
and for sufficiently small $h$ introduce  \smallskip
$$
w(t,x,\xi)=v(\theta^{0}+h,t,x+\xi)-v(\theta^{0},t,x),
\smallskip$$
where $t,x,\xi\in \bar Q$ with  
$$
Q:=[0,1)
\times\big\{(x,\xi):x,x+\xi\in B_{1}\big\}.
\medskip$$
Note for the future that, by embedding theorems
$D_{x}v(\theta,t,x)$ is a continuous function in $\bar C_{1} $
for any $\theta$.

Next, observe that   \smallskip
$$
F\big(\theta^{0}+h,D^{2}_{x}v(\theta^{0}+h,t,x+\xi),t\big)
-F\big(\theta^{0},D^{2}_{x}v(\theta^{0},t,x),t\big)
\medskip$$
$$
=\big[
F\big(\theta^{0},D_{x}^{2}v(\theta^{0}+h,t,x+\xi),t\big)
-F(\theta^{0},D^{2}_{x}v(\theta^{0},t,x),t\big)\big]+I
\medskip$$
$$
=a^{ij}D _{x^{i}x^{j}}w+I,
\medskip$$
 in $Q$ (a.e.)
where $(a^{ij})$ is a measurable $\bS_{\delta}$-valued function and  \smallskip
$$
I=F\big(\theta^{0}+h,D^{2}_{x}v(\theta^{0}+h,t,x+\xi),t\big)
-F\big(\theta^{0},D_{x}^{2}v(\theta^{0}+h,t,x+\xi),t\big).
\medskip$$
Since $|D_{\theta}F|\leq \chi+|\sfu''|$ by assumption, we have  \smallskip
$$
I=\tau h+h \varepsilon^{ij}D_{x^{i}x^{j}}
v(\theta^{0}+h,t,x+\xi)
=\tau h+h \varepsilon^{ij}D_{x^{i}\xi^{j}}
w( t,x,\xi),
\medskip$$
where $|\tau|\leq \chi$ and $(\varepsilon^{ij})$ is an $\bS$-valued function
with norm majorated by one. Hence,  in $Q$ (a.e.)
 we have \smallskip
$$
\partial_{t}w+\cL w+\tau h=0,
$$
where
$$
\cL=\big[a^{ij}-N_{0}h^{2}\delta^{ij}\big]D^{2}_{x^{i}x^{j}} +
h \varepsilon^{ij}D_{x^{i}\xi^{j}} 
+N_{0}h^{2}\delta^{ij}D_{\xi^{i}\xi^{j}} 
\medskip$$
is a uniformly elliptic operator for an appropriate $N_{0}=N_{0}(\delta,d)$
and all sufficiently small $h\ne0$.

Notice that on $\partial'Q$ either $t=1$, and then  \smallskip
$$
|w|\leq \sup_{|x|\leq 1}|Dg|\,|\xi|=\sup_{|x|\leq 1} \big|Dg-Dg(1,0) \big|\,|\xi|,
$$
  or $|x|\leq 1$ and
$|x+\xi|=1$, in which case $v(\theta^{0}+h,t,x+\xi)=
v(\theta^{0} ,t,x+\xi) $ and
$$
\big|w(t,x,\xi)\big|\leq \sup_{  C_{1}}\big|D_{x}v(\theta^{0},\cdot)\big|\,|\xi|,
$$
or else $|x+\xi|\leq 1$ and
$|x |=1$, in which case $v(\theta^{0} ,t,x )=
v(\theta^{0}+h ,t,x ) $ and   \smallskip
$$
\big|w(t,x,\xi)\big|\leq \sup_{  C_{1}}\big|D_{x}v(\theta^{0}+h,\cdot)\big|\,|\xi|.
$$
In all cases on $\partial' Q$ we have \smallskip
$$
\big|w(t,x,\xi)\big|\leq N_{1}|\xi|\leq N_{1}h +N_{1}
 h ^{-1}|\xi|^{2},
$$
where
$$
N_{1}= [Dg]_{C^{\kappa}(C_{1})}+\max_{ \theta=\theta^{0},\theta^{0}+h}
\sup_{  C_{1}}\big|D_{x}v(\theta ,\cdot)\big|.
$$

As is easy to see, there is a constant $N_{2}=N_{2}(d,\delta)$
such that, for the function
$$
\phi(t,x,\xi)=h\chi(1-t)+N_{1} h 
+N_{1}  h ^{-1}\big[|\xi|^{2}+N_{2}h^{2}\big(1-|x|^{2}\big)\big]
$$
we have
$$
\partial_{t}\phi+\cL\phi+\tau h\leq0
\medskip$$
in $Q$ and, of course, $|w|\leq\phi$ on $\partial'Q$.
By the maximum principle   $|w|\leq\phi$
in $\bar Q$, in particular, (take $\xi=0$) in $\bar C_{1}$  \smallskip
$$
\big|v(\theta^{0}+h,t,x )-v(\theta^{0},t,x)\big|
\leq h \chi +NN_{1}h,
\medskip$$
where $N=N(\delta,d)$.

It follows that to prove the lemma, it suffices to show that
for any $\theta\in[-\theta_{0},\theta_{0}]$,
with a constant $N=N(\delta,d,\kappa)$, we have  \smallskip
\begin{equation}
                                                \label{12,18.5}
\sup_{  C_{1}}\big|D_{x}v(\theta ,\cdot)\big|
\leq N [Dg]_{C^{\kappa}(C_{1})}+N \sup_{t}\big|F(\theta,0,t)\big|.
\medskip\end{equation}

By applying finite-difference operators with respect to $x$
to \eqref{12,18.1}, we see that, for small $h$
and unit $l\in\bR^{d}$, the function $ \big[v(\theta,t,x+hl)
-v(\theta,t,x)\big]/h$
satisfies a parabolic equation with zero free term in a domain
slightly smaller than $C_{1}$. Hence, its sup over the domain
is achieved on the parabolic boundary. By letting $h\to0$
we conclude that  \smallskip
$$
\sup_{  C_{1}}\big|D_{x}v(\theta ,\cdot)\big|=
\sup_{\partial' C_{1}}\big|D_{x}v(\theta ,\cdot)\big|,
\medskip$$
and since $v(\theta,1,x)=g(1,x)$ for $|x|\leq1$,
to prove \eqref{12,18.5}, it suffices to prove that
\begin{equation}
                                                \label{12,27.3}
\sup_{(0,1)\times\partial B_{1}}\big|D_{x}v(\theta ,\cdot)\big|
\leq N [Dg]_{C^{\kappa}(C_{1})}+N \sup_{t}\big|F(\theta,0,t)\big|.
\medskip\end{equation}

We fix $\theta$ and observe that  \smallskip
$$
0=\partial_{t}v(\theta,t,x)+\big[F\big(\theta,D^{2}v(\theta,t,x),t\big)
-F(\theta,0,t)\big]+f 
\medskip$$
$$
=\partial_{t}v(\theta,t,x)+a^{ij}D_{ij}v(\theta,t,x)+f ,
\medskip$$
(a.e.), where $(a^{ij})$ is a measurable $\bS_{\delta}$-valued
function and $f =F(\theta, 0,t)$. After that
\eqref{12,27.3} immediately follows from
Lemma \ref{lemma 1.1.5}.   The lemma is proved.     \qed
                                              
 {\bf Proof of Theorem \ref{theorem 12,28.1}}. For $\theta\in
[-\theta_{0},\theta_{0}]$ introduce  \smallskip
$$
F(\theta,\sfu'',t)=\max\Big(F_{0}(\sfu'',t)+
\theta\big(K_{1}/\theta_{0}+|\sfu''|\big),P(\sfu'')-K\Big).
\medskip$$
By Theorem \ref{theorem 4,5,1}
  equation \eqref{12,18.1} in $C_{R}$   
with boundary condition $v=g$ on $\partial' C_{R}$
admits a solution 
$v=v(\theta,\cdot)\in W^{1,2}_{p}(
C_{R}) $ for any $p\geq1$. By the maximum principle the solution is
unique. Obviously, $\big|D_{\theta}F (\theta,\sfu'')\big|\leq
K_{1}/\theta_{0}+|\sfu''|$ whenever the left-hand side is well defined.
Also $\big|F(\theta,0,t)\big|\leq K_{1}$. After that,
it only remains to observe that,
 in light of Lemma \ref{lemma 12,18.1},   \smallskip
$$
\big|v^{(+)}-v^{(-)}\big|\leq
\int_{-\theta_{0}}^{\theta_{0}}
\Big(R^{2} \big(K_{1}/\theta_{0} +NK_{1} \big)  +NR^{\kappa}
[g]_{C^{\kappa}(C_{r})}\Big)\,d\theta
\medskip$$
$$
\leq NR^{2} K_{1}  +NR^{\kappa}\theta_{0}
[g]_{C^{\kappa}(C_{r})}.
\medskip$$

The theorem is proved.  \qed

Later on we will need one more piece
of information about $v^{(\pm)}$, before which
we prove the following two auxiliary facts.
 
\begin{lemma}
                                      \label{lemma 12,18.2}
Take $k\in\{1,...,d\}$, $\rho\in\bR$, $\lambda\in(0,\infty)$, and introduce \smallskip
$$
V (t,x)=\exp\big(-\lambda|x^{k}-\rho|^{2}/(r^{2}-t)\big).
\medskip$$
Let $(a^{ij})\in\bS$ be such that   $0\leq a^{kk}
\leq 1/(4\lambda)$. Then for $t< r^{2}$ it holds that\smallskip
\begin{equation}
                                           \label{12,18.4}
\partial_{t}V(t,x)+a^{ij}D_{ij}V(t,x)\leq0.
\medskip\end{equation}
\end{lemma}

The proof is achieved by a straightforward
computation showing that the left-hand side of
\eqref{12,18.4} equals   \smallskip
$$
V (t,x)\Big[-\lambda\frac{|x^{k}-\rho|^{2}}{(r^{2}-t)^{2}}
+4\lambda^{2} a^{kk}\frac{|x^{k}-\rho|^{2}}{(r^{2}-t)^{2}}
-2\lambda a^{kk}\frac{1}{r^{2}-t}\Big]\leq0.
\medskip$$

\begin{lemma}
                                         \label{lemma 12,9.1}
Let $r\in(0,\infty)$, $v\in W^{1,2}_{d+1,\loc}(C_{r})\cap C(\bar C_{r})$,
and assume that, for   constants $\gamma\in(0,2]$, $M\geq0$, we have  \smallskip
$$
\big|v(r^{2},x)\big|\leq M|x|^{\gamma}\quad\text{if}\quad |x|\leq r.
\medskip$$
Also assume that $\big|\partial_{t}v+a^{ij}D_{ij}v\big|\leq\theta$ in $C_{r}$
for an $\bS_{\delta}$-valued function $a=(a^{ij})$ and a constant $\theta\in[0,\infty)$. Then there exists
a constant $N=N(d,\delta,\gamma)$ such that for $t\in[0,r^{2}]$  \smallskip
\begin{equation}
                                                         \label{12,9.1}
\big|v(t,0)\big|\leq N\big(M+r^{-\gamma}\sup_{C_{r}}
|v|\big) (r^{2}-t )^{\gamma/2}
+\theta (r^{2}-t ).
\medskip\end{equation}
\end{lemma}

Proof. Take  $\varepsilon >0$, set $\rho=r/\sqrt{d}$, 
$\lambda=\delta/4$,
  and consider the function   \smallskip
$$
V(t,x)=M\varepsilon+M\varepsilon^{1-2/\gamma}\big(|x|^{2}
+2(r^{2}-t)d/\delta\big)
+\theta(r^{2}-t)
\medskip$$
$$
+\sup_{C_{r}}|v|
\sum_{k=1}^{d} \big(V_{(+)}^{k}(t,x)+V_{(-)}^{k}(t,x) \big) ,
\medskip$$
where
$$
V_{(\pm)}^{k}(t,x)=\exp\big(-\lambda|x^{k}\pm\rho|^{2}/(r^{2}-t)\big).
\medskip$$
Observe that,  by Lemma \ref{lemma 12,18.2},
we have  \smallskip
$$
\partial_{t}V  (t,x)+a^{ij}D_{ij}V\leq-\theta  
\medskip$$
in $C_{r}$, in particular, in the cylinder $(0,r^{2})\times(-\rho,\rho)^{d}$.
On the parabolic boundary of this cylinder either $t=r^{2}$
and $V\geq M|x|^{\gamma}\geq |v|$ by Young's inequality
($\gamma\leq 2$), or
one of $x^{i}$ is equal to either $\rho$ or $-\rho$, when the corresponding
$V^{i}$ equals 1, and again $V\geq |v|$. By the maximum principle,
for $t\in[0,r^{2}]$ we have \smallskip
$$
 v(t,0)\leq V(t,0)\leq NM\big(\varepsilon+ \varepsilon^{1-2/\gamma} (r^{2}-t)\big)
\medskip$$
$$
+\theta(r^{2}-t)
+N\sup_{C_{r}}|v|\exp\big(-\lambda\rho^{2}/(r^{2}-t)\big).
\medskip$$
After that, to estimate $v(t,0)$ by the right-hand side of \eqref{12,9.1}, it only remains to take the inf with respect to 
$\varepsilon>0$ and observe that  \smallskip
$$
\exp\big(-\lambda\rho^{2}/(r^{2}-t)\big)\leq N\big(r^{2}-t\big)^{\gamma/2}/
r^{\gamma}.
\medskip$$
Similarly $-v(t,0)$ is estimated from above.   
The lemma is proved.   \qed

Recall that $v^{(\pm)}$ are introduced before Theorem \ref{theorem 12,28.1}.

\begin{theorem}
                                         \label{theorem 3,15,1}   
There exist constants $\kappa_{0}
=\kappa_{0}(d,\delta)\in(1,2 ) $ and $N\in(0,\infty)$
depending only on $d$ and $\delta$
such that for any $r\in(0,R)$  \smallskip
$$
\big[v^{(\pm)}\big]_{C^{ \kappa_{0} }(  C_{r})}
\leq N(R-r)^{-\kappa_{0}}\big[\osc_{\partial' C_{R}}(g-\hat{g})
+K_{1}R^{2 }\big],
\medskip$$
where $\hat{g}=\hat{g}(x)$ is any affine function of $x$.
\end{theorem}

Proof. As usual we may take $\hat g=0$ and recall that
$Dv^{(\pm)}$ is bounded and even H\"older continuous in $C_{R}$
since $v^{(\pm)}\in W^{1,2}_{p}(C_{R})$ for any $p\geq1$.
Then  
  observe that for any $\gamma\in(0,1)$ and
function $f(x)$ of one variable
$x\in[0,\varepsilon]$, $\varepsilon>0$,  we have  \medskip
$$
\big|f'(0)\big|\leq \big|f'(0)-\big(f(\varepsilon)-f(0)\big)/\varepsilon\big|+
\varepsilon^{-1}\osc_{[0,\varepsilon]}f
\leq\varepsilon^{\gamma}[f']_{C^{\gamma}[0,\varepsilon]}+
\varepsilon^{-1}\osc_{[0,\varepsilon]}f.
\medskip$$
By applying this fact to functions $v(x)$ given in $B_{R}$ 
we obtain that for any $r_{n+1}<r_{n+2} \leq R$
and  any $\varepsilon\in(0,1)$ and smooth $v=v(x)$  \smallskip
\begin{equation}
                                              \label{11.1.1}
|Dv |\leq\varepsilon^{\gamma}(r_{n+2}-r_{n+1 })^{\gamma} 
[Dv ]_{C^{\gamma}(  B_{r_{n+2}})}
+\varepsilon^{-1}(r_{n+2}-r_{n+1})^{-1}
 \osc_{  B_{R}}v 
\end{equation}
in $\bar B_{r_{n+1}}$.

Next, for
 unit $l\in\bR^{d}$ and $h>0$ define  \smallskip
$$
\delta_{h,l}  u(t,x)=  
h^{-1}\big[u(t,x+hl)-u(t,x)\big]
\medskip$$
 and note that
for any $ r_{1}\in(0,R) $ the function
$\delta_{ h,l }v^{(+)}$ satisfies an equation of the type  \smallskip
$$
\partial_{t}\delta_{h,l}v^{(+)}+a^{ij}D_{ij}\delta_{h,l}v^{(+)}=0
\medskip$$
in $C_{r_{1}}$ (a.s.)
with some measurable $(a^{ij})$ taking values in $\bS_{\bar\delta}$
 if $h$ is sufficiently small. By the Krylov-Safonov theorem,
 for $r_{0}\in(0,r_{1})$ and
perhaps even smaller  $h$   we have  \smallskip
$$
\big[\delta_{h,l}v^{(+)}\big]_{C^{ \gamma}(  C_{ r_{0}})}
\leq N(r_{1}-r_{0})^{-\gamma}\sup_{  C_{ r_{1} }}\big|\delta_{h,l}v^{(+)}\big|,
\smallskip$$
where $\gamma\in(0,1)$ and $N$ depend only on $\delta$
and $d$. By letting $h\to0$ we conclude   \smallskip
\begin{equation}
                                               \label{9.4.1}    
 [Dv^{(+)} ]_{C^{ \gamma}(  C_{ r_{0}})}
\leq N(r_{1}-r_{0})^{-\gamma}\sup_{  C_{ r_{1} }} |Dv^{(+)} |.
\end{equation}

By using \eqref{11.1.1} and \eqref{9.4.1} and setting
$$
r_{0}=r,\quad
r_{n}=r+(R-r)\sum_{k=1}^{n}2^{-k},\quad n\geq1,
$$
we   obtain  \smallskip
$$
A_{n}:=\sup_{[0,r_{n}^{2}]}\big[Dv^{(+)}(t,\cdot)
\big]_{C^{ \gamma}(  B_{ r_{n}})}
\leq N(r_{n+1}-r_{n})^{-\gamma}
\sup_{  C_{ r_{n+1}}}|Dv^{(+)}|
\medskip$$
\begin{equation}
                                               \label{9.4.02}
\leq N_{1}\varepsilon^{\gamma}  
A_{n+2}+N_{2}(R-r)^{-(1+\gamma)}\varepsilon^{-1}2^{(1+\gamma)n}
\osc_{  C_{R}}v^{(+)},
\medskip\end{equation}
where the constants $N_{i}$ are different from the one in \eqref{9.4.1}
but still depend  only on $\delta$ and $d$.
Without losing generality we may assume that $N_{1}\geq1$
and we first take $\varepsilon$ so that  \smallskip
$$
N_{1}\varepsilon^{\gamma} 
=2^{-5},
\smallskip$$
then take $n=2k$, $k=0,1,...$,  multiply both parts
of \eqref{9.4.02} by $2^{-5k}$ and sum up with respect to $k$.
Then upon observing that $(1+\gamma)2k\leq 4k$ we get  \smallskip
$$
\sum_{k=0}^{\infty}A_{2k}2^{-5 
k}\leq
\sum_{k=1}^{\infty}A_{2k}2^{-5k}+N(R-r)^{-(1+\gamma)}
\sum_{k=0}^{\infty}
2^{-k}\osc_{  C_{R}}v^{(+)}.
\medskip$$
By canceling (finite) like terms we find    \smallskip
\begin{equation}
                                               \label{9.4.03}
\sup_{[0,r^{2}]}\big[Dv^{(+)}(t,\cdot)\big]_{C^{ \gamma}( B_{r})}
\leq N(R-r)^{-(1+\gamma)}\osc_{ C_{R}}v^{(+)}.
\medskip\end{equation}

Note for the future that 
  \eqref{9.4.03} and the second inequality in \eqref{9.4.02}
also imply that
\begin{equation}
                                         \label{12,11.01} 
 \sup_{ C_{ r }} |Dv^{(+)} |\leq
N(R-r)^{-1}\osc_{  C_{ R}}v^{(+)}.
\medskip\end{equation}

Next, we use the fact that $v^{(+)}$ itself satisfies the equation \vspace{5pt}
\begin{equation}
                                                    \label{1,3.1}
0=\partial_{t}v^{(+)}+\max\big(F^{(+)} [v^{(+)} ],
P [v^{(+)} ]-K\big)
-\max(K_{1},-K)+K_{1}
 \end{equation}  \vspace{0pt}
$$
=\partial_{t}v^{(+)}+a_{ij}D_{ij}v^{(+)} 
+K_{1}
\medskip$$
in $C_{R}$ (a.e.) with some measurable $(a_{ij})$ taking values in $\bS_{\bar\delta}$.
 Furthermore, for any $(t_{0},x_{0})\in C_{R}$
the function   \smallskip
$$
v(t,x):=v^{(+)}(t,x)-v^{(+)}(t_{0},x_{0})-
(x^{i}-x^{i}_{0 })D_{i}v^{(+)}(t_{0},x_{0})
\medskip$$
satisfies the same equation and, owing to
\eqref{9.4.03},   \vspace{5pt}
$$
|v(t_{0},x) |
\leq \big[Dv^{(+)}(t_{0},\cdot)\big]_{C^{\gamma}
(  B_{r})}
|x-x_{0}|^{1+\gamma}
\leq N(R-r)^{-(1+\gamma)}|x-x_{0}|^{1+\gamma}
\osc_{  C_{R}}v^{(+)}
\medskip$$
if $(t_{0},x_{0})\in C_{r}$
and $|x-x_{0}|\leq \rho:=(R-r)/2$.

Also, for $t_{1}=0\vee \big(t_{0}-\rho^{2}\big)$,
due to \eqref{12,11.01}, we have that  \vspace{5pt}
$$
\sup_{C_{t_{0}-t_{1},\rho}(t_{1},x_{0})} |v| 
\leq \osc_{C_{t_{0}-t_{1},\rho}(t_{1},x_{0})}v^{(+)}
+N\rho\big|D v^{(+)}(t_{0},x_{0})\big|\leq N
\osc_{  C_{ R}}v^{(+)}.
 $$  \vspace{5pt}

Then we apply Lemma \ref{lemma 12,9.1} with
$C_{t_{0}-t_{1},\rho}(t_{1},x_{0})$ in place of $C_{r}$
and take into account that $\rho^{- (1+\gamma) }\leq N(R-r)^{-(1+\gamma)}$.  Then
for $t\in[t_{1},t_{0}]$ we obtain    \vspace{5pt}
$$
\frac{\big|v^{(+)}(t_{0},x_{0})-v^{(+)}(t,x_{0})\big|}
{(t_{0}-t)^{(1+\gamma)/2}}\leq
\frac{N}{(R-r)^{ 1+\gamma }}
\osc_{ C_{R}}v^{(+)}+K_{1}(t_{0}-t)^{(1-\gamma)/2}.
\medskip$$   \vspace{0pt}

Here $t_{0}-t\leq \rho^{2}$ and $(t_{0}-t)^{(1-\gamma)/2}\leq
R^{2}(R-r)^{-(1+\gamma)}$. Therefore  \vspace{5pt}
\begin{equation}
                                         \label{12,11.2}
\frac{\big|v^{(+)}(t_{0},x_{0})-v^{(+)}(t,x_{0})\big|}
{(t_{0}-t)^{(1+\gamma)/2}}\leq
\frac{N}{(R-r)^{ 1+\gamma }}\big[
\osc_{  C_{R}}v^{(+)}+K_{1}R^{2}\big].
 \end{equation}\vspace{5pt}

If $x_{0}$, $t_{0}$ are as above but $0\leq t<t_{1}$,
then $t_{0}-t\geq \rho^{2}=(R-r)^{2}/4$ and  \vspace{5pt}
$$
\big|v^{(+)}(t_{0},x_{0})-v^{(+)}(t,x_{0})\big|
\leq \osc_{ C_{R}}v^{(+)}
= (R-r)^{-(1+\gamma)}\big(4\rho^{2}\big)^{(1+\gamma)/2}
\osc_{  C_{ R}}v^{(+)},
\medskip$$
so that \eqref{12,11.2} holds again.

This provides the necessary estimate of the oscillation of $v^{(+)}$
in the time variable and along with \eqref{9.4.03} shows that  \smallskip
$$
 [v^{(+)} ]_{C^{ 1+\gamma}( C_{r})}
\leq N(R-r)^{-(1+\gamma)}\big[\osc_{  C_{ R}}v^{(+)}
+K_{1}R^{2}\big].
$$

Now the assertion of the theorem about $v^{(+)}$
with $\kappa_{0}=1+\gamma$ follows
from the fact that  \smallskip
$$
\osc_{  C_{R}}v^{(+)}\leq \osc_{\partial' C_{R}}g+NK_{1}R^{2}.
$$

The function $v^{(-)}$ is considered similarly
with the only difference that it satisfies an equation
similar to \eqref{1,3.1} with $-\max(-K_{1},-K) +\max(-K_{1},-K)$
in place of $-\max( K_{1},-K)+K_{1}$. Since 
$|\max(-K_{1},-K)|\leq K_{1}$ this difference
is irrelevant. The theorem is proved.\qed

\mysection
{Main estimate for solutions of an auxiliary    
cut-off equation}
                                               \label{section 011.4.4}
Let $F\big(\sfu'_{0},\sfu'',t,x\big)$
satisfy Assumptions \ref{assumption 12,29.2}     
(i), (iii), (iv),   and (v).
Take $K\in(0,\infty)$, $R\in (0,R_{0}]$,
 and $g\in W^{1,2}_{\infty}(C_{R}) $ and  take $P[u]$ as in the beginning 
of Section~\ref{section 12,31.1}.
 
By Theorem \ref{theorem 4,5,1}  there exists
  $u\in W^{1,2}_{p}(C_{R})$ for all $p\geq 1$,
such that $u=g$ on $\partial'C_{R}$ and the equation 
\begin{equation}
                                          \label{12,31.1}
\partial_{t}u+\max\big(F[u],P[u]-K\big)=0.
\end{equation}
holds (a.e.) in $C_{R}$.   

Here is the  result of this section.
\begin{theorem}
                                        \label{theorem 12,31.1}
Suppose that 
$$
\osc_{\partial'C_{R}} g < \hat R_{0}.
$$
Then there exist a constant  $\kappa_{0}=\kappa_{0}(d,\delta)
\in(1,2]$ such that, if $r<R\leq R_{0}$,
one can find an affine function $\hat{u}=\hat{u}(x)$
for which
\begin{align*}
|u-\hat{u}|\leq &\, N K_{1}R^{2}+
N\theta_{0}[g]_{C^{ \kappa}( C_{R})}R^{\kappa}
\\[10pt]
 &\, +
Nr^{ \kappa_{0}}(R-r)^{-\kappa_{0}}
\big[\osc_{\partial' C_{R}}(g-\hat{g})
+K_{1}R^{2 }\big]
\end{align*}

\noindent
  in $\bar C_{r}$ for any $\kappa\in(0,2)$, where
  $\hat{g}=\hat{g}(x)$ is any affine function of $x$ 
  and $N=N(d,\delta,\kappa)$. 

\end{theorem}

Proof.  Set 
$$
F_{0}(\sfu'',t)=F\big(g(R^{2},0),\sfu'',t,0\big) ,
\medskip$$
 and take $v^{(\pm)}$ 
from Section \ref{section 12,31.1}. 
 
Observe that since $\max\big(F\big(\sfu'_{0},0,t,x\big),P(0)-K\big)=0$,
as we have seen many times in the past,    \smallskip
$$
\partial_{t}u+a^{ij}D_{ij}u=0
\medskip$$
for an appropriate $\bS_{\bar\delta}$-valued $(a^{ij})$.
It follows that    \smallskip
$$
\osc_{ C_{R}} u =\osc_{\partial'C_{R}} g< \hat R_{0},
\medskip$$
and in $C_{R}$ by Assumption \ref{assumption 12,29.2} (iv)  \smallskip
$$
F\big(u(t,x),D^{2}u(t,x),t,x\big)\leq F_{0}\big(D^{2}u(t,x),t\big)+
K_{1}+\theta_{0}\big|D^{2}u(t,x)\big|,
\medskip$$
$$
 F_{0}\big(D^{2}u(t,x),t\big)-K_{1}-\theta_{0}\big|D^{2}u(t,x)\big| \leq 
F\big(u(t,x),D^{2}u(t,x),t,x\big).
\medskip$$
Since the operators $F^{(\pm)}$ satisfy the maximum principle,
the above inequalities imply that $v^{(-)}\leq u\leq v^{(+)}$.

Next, set  $p^{(\pm)}(x)=v^{(\pm)}(0)
+x^{i}D_{i}v^{(\pm)}(0)$
and observe that by Theorem~\ref{theorem 3,15,1}  in $\bar C_{r}$
we have \smallskip
$$
|p^{(+)}-v^{(+)}|+|p^{(-)}-v^{(-)}|\leq
Nr^{\kappa_{0}}
 N(R-r)^{-\kappa_{0}}\big [\osc_{\partial' C_{R}}(g-\hat{g})
+K_{1}R^{2 }\big] ,
\medskip$$
  where $\kappa_{0}$ is
taken from Theorem \ref{theorem 3,15,1}. Furthermore
by Theorem \ref{theorem 12,28.1} in $\bar C_{R}$  \smallskip
$$
|p^{(+)}-p^{(-)}|\leq|p^{(+)}-v^{(+)}|+|p^{(-)}-v^{(-)}|
\medskip$$
$$
+NR^{2} K_{1} +NR^{\kappa}\theta_{0} 
[g]_{C^{\kappa}(C_{R})}.
\medskip$$

After that it only remains to note that, since
$v^{(-)}\leq u\leq v^{(+)}$, we have  \smallskip
$$
p^{(+)}-v^{(+)}\leq p^{(+)}-u\leq 
 [p^{(+)}-p^{(-)}]+[p^{(-)}-v^{(-)}]
+v^{(-)}-u
\medskip$$
$$
\leq [p^{(+)}-p^{(-)}]+[p^{(-)}-v^{(-)}].
\medskip$$
 
The theorem is proved.    \qed

\mysection
{Proof of Theorem
\protect\ref{theorem 12,30.1}}
                                      \label{section 1,3.1}

We first    take
a constant $K>0$, assume that $g\in W^{1,2}_{\infty}(\bR^{d+1})$,
$\bar G$ is bounded,
   and investigate
solutions $v_{K}$ of the cut-off equation
\eqref{9.28.4}.
Take $\kappa_{0}=\kappa_{0}(d,\delta)\in(1,2]$ from Theorem   
\ref{theorem 12,31.1}.    
Naturally, we suppose that the assumptions
of Theorem
 \ref{theorem 12,30.1}:
Assumptions \ref{assumption 12,29.1},
\ref{assumption 12,29.2}, and \ref{assumption 12,30.3},
are satisfied with $\theta_{0}$ yet to be specified.

\begin{lemma}
                                       \label{lemma 1,3.2}
Let $R\in(0,R_{0}]$  and let
$ v \in W^{1,2}_{\infty}( C_{R}) $ be a   solution
of \eqref{9.28.4} in $\bar C_{R}$ such that
\begin{equation}
                                           \label{1,24,4}
\osc_{\partial' C_{R}} v< \hat R_{0}.
\medskip\end{equation}
Then for each $r\in(0,R)$
  one can find
 an affine function $\hat{v}(x)$
such that in $C_{r}$ for any $\kappa\in(1,2)$  
$$
|v-\hat{v}|\leq  
N\theta_{0}[v]_{C^{ \kappa}(  C_{R})}R^{\kappa}+
Nr^{ \kappa_{0}}(R-r)^{-\kappa_{0}}  
\big[R^{\kappa}[v]_{C^{ \kappa}(  C_{R})}+K_{1}R^{2}\big]
\medskip$$
$$
 +
NK_{0}R^{2}\sup_{ C_{R}}\big (|v|+|Dv|\big ) 
  +N  K_{1} R^{2} +NR^{2-(d+2)/p}
 \|\bar G\|_{L_{p}(C_{R})},  
\medskip$$
where   the
constants $N$ depend  only on $d$, $\delta$,  and $\kappa$.

\end{lemma}

Proof. 
Observe that  \smallskip
$$
\max\big (H[v] ,P[v]-K\big )=\max\big( F[v] ,P[v]-K\big )+h
\medskip$$
where $h$, defined by the above equality, satisfies  \smallskip
$$
|h|\leq\big|H[v]-F[v]\big| \leq K_{0}\big(|v|+|Dv|\big)+\bar{G} .
\medskip$$
Next define $ u \in W^{1,2}_{d+1}(C_{R})$
 as a   solution  of  \smallskip   
$$
\partial_{t}u+\max\big(F[u],P[u]-K\big)=0
\medskip$$
with  boundary data $u=v$ on $\partial' C_{R}$,
which exists by Theorem \ref{theorem 4,5,1}.
Then,  in light of the fact that  \smallskip
$$
\max\big(F(u(t,x),0,t,x),P(0)-K\big)=0,
\medskip$$
there exists an 
$\bS_{\bar\delta}$-valued measurable function $a$ such that
in $C_{R}$ (a.e.) we have    \smallskip
$$
\partial_{t}(v-u)+a^{ij}D_{ij}(v-u)+h=0.
\medskip$$
By the parabolic Aleksandrov estimate  \smallskip
$$
|v-u|\leq NR^{d/{d+1}}\|h\|_{L_{d+1}(C_{R})}=NR^{2}
\Big(\dashint_{C_{R}}|h|^{d+1}\,dxdt\Big)^{1/(d+1)} 
\medskip$$
$$
\leq NK_{0}R^{2}\sup_{ C_{R}}\big(|v|+|Dv|\big)+NR^{2-(d+2)/p}
 \|\bar G\|_{L_{p}(C_{R})}.
\medskip$$

After that our assertion follows from 
  Theorem \ref{theorem 12,31.1} and the lemma is proved.  
\qed

We need a
\index{characterization of $C^{1+\alpha}$-functions}%
characterization of $C^{1+\alpha}$-functions.
\begin{lemma} 
                                      \label{lemma 10.23.1}
Let $r_{0}\in(0,\infty)$, $\kappa\in(1,2)$, $\phi\in 
C^{\kappa}(  C_{r_{0}})$ and assume that there is a constant
$N_{0}$ such that for any $(t,x)\in C_{r_{0}}$ 
and $r\in(0,2r_{0}]$
there exists an affine function $\hat{\phi}=\hat{\phi}(x)$ such that
$$
\sup_{  C_{r}(t,x)\cap
  C_{r_{0}}} |\phi-\hat{\phi}|\leq N_{0}r^{\kappa}.
$$
Then
$$
[\phi]_{C^{\kappa}(  C_{r_{0}})}\leq NN_{0},
\medskip$$
where $N  $ depends only on $d$ and $\kappa$.

\end{lemma}

Proof. The fact that, for any $t\in(0,r^{2})$, we have  \smallskip
$$
\big[D\phi(t,\cdot)\big]_{C^{\kappa-1}(  B_{r_{0}})}\leq NN_{0}
\medskip$$
follows from Theorem 2.1 of \cite{Sa_88}. To estimate
$|\phi(t,x)-\phi(s,x)|$ we may assume that $t>s$, so that
$(t,x),(s,x)\in \bar C_{r}(s,x)$, where $r=\sqrt{t-s}$.
Then, for an appropriate $\hat\phi(x)$  \smallskip
$$
\big|\phi(t,x)-\phi(s,x)\big|\leq \big|\phi(t,x)-\hat\phi(x)\big|
+\big|\phi(t,x)-\hat\phi(x)\big|\leq 2N_{0}r^{\kappa}=2N_{0}(t-s)^{\kappa/2}.
\medskip$$
The lemma is proved. \qed

 \begin{lemma} 
                                          \label{lemma 1,3.4}
Take $r_{1}\in(0,R_{0}]$, $r_{0}\in(0,r_{1})$,
and define  \smallskip
$$
\kappa  =\kappa(d,\delta,p)  =
\frac{1+\kappa_{0}}{2}\wedge\Big(2-\frac{d+2}{p}\Big).
\medskip$$
Let 
$ v \in W^{1,2}_{\infty}(C_{r_{1}}) $ be a 
  solution
of \eqref{9.28.4} in $C_{r_{1}}$.
 Then there exists $\theta_{0}=\theta_{0}( d,\delta)\in(0,1]$ such that,
if Assumption \ref{assumption 12,29.2} \(iv\,\)
 is satisfied with this
$\theta_{0}$ and  \smallskip
\begin{equation}
                                           \label{01,24,4}
\osc_{  C_{r_{1}}} v< \hat R_{0},
\end{equation}
  then     
\begin{align}
[v]_{C^{\kappa}(  C_{r_{0}})}
\leq &\, (1/2)[v]_{C^{  \kappa}( C_{r_{1}})}
+N(K_{0}+1) (r_{1}-r_{0})^{-\kappa }\sup_{  C_{r_{1}}} |v|
\nonumber\\[10pt]
                                               \label{1,3.5}
 &\, +N (K_{0}+1)(r_{1}-r_{0})^{-(\kappa-1) }
\sup_{  C_{r_{1}}} |Dv| 
+  N ( K_{1}+\|  G\|_{L_{p}(C_{r _{1}})} )  ,
\end{align}
where $N =N (d,\delta )$. 
\end{lemma}

Proof. 
Take $(t_{0},x_{0})\in C_{r_{0}}$, $\varepsilon\in(0,1)$
to be specified later, define   \smallskip
$$
r'_{0}=\frac{\varepsilon}{3}(r_{1}-r_{0}), 
\medskip$$
and notice that for any $(t,x)\in C_{r_{0}'}(t_{0},x_{0})$,
  $r\in (0,2r_{0}']$, and $R=\varepsilon^{-1}r$,   we have
$R\leq r_{1}\leq R_{0}\leq 1$ and
$$
C_{R}(t,x)\subset C_{r _{1}} .
\medskip$$
  Therefore,
by Lemma \ref{lemma 1,3.2} we can find an affine function
$\hat{v}(x)$ such that   \smallskip
$$
\sup_{C_{r}(t,x)\cap C_{r'_{0}}(t_{0},x_{0})}
|v-\hat{v}|\leq
\sup_{  C_{r}(t,x) }
|v-\hat{v}|\leq N\theta_{0}[v]_{C^{ \kappa}(  C_{R}(t,x))}\varepsilon ^{-\kappa }
r^{\kappa }
\vspace{5pt}$$
$$
+
N\varepsilon^{ \kappa_{0}-\kappa }
(1-\varepsilon)^{-\kappa_{0}}
r^{\kappa }[v]_{C^{ \kappa}(  C_{R}(t,x))}
+N\varepsilon^{\kappa_{0}-2}(1-\varepsilon)^{-\kappa_{0}}K_{1}r^{2}
\vspace{15pt}$$
$$
 +
NK_{0}\varepsilon^{-2}
r^{2}\sup_{ C_{R}(t,x)}\big(|v|+|Dv|\big) +N  K_{1} 
\varepsilon^{-2}r^{2}
\vspace{10pt}$$
$$
+N\varepsilon^{-\kappa}r^{\kappa}
\|\bar G\|_{L_{p}(C_{r _{1}})}
 \leq N r^{\kappa }I(\theta_{0},\varepsilon,r_{1}),  
\medskip$$

\noindent
where the constants $N$ depend only on $d$  and $\delta$ and \smallskip
$$
I(\theta_{0},\varepsilon,r_{1}):=\big(\theta_{0}
\varepsilon ^{-\kappa }+\varepsilon^{ \kappa_{0}-\kappa }
(1-\varepsilon)^{-\kappa_{0}}
\big)[v]_{C^{  \kappa}(  C_{r_{1}})}  
+\varepsilon^{\kappa_{0}-2}(1-\varepsilon)^{-\kappa_{0}}K_{1}
 $$  \vspace{0pt}
$$
+\varepsilon^{-2}K_{0}\sup_{  C_{r_{1}}}\big(|v|+|Dv|\big)
+ K_{1} 
\varepsilon^{-2}+\varepsilon^{-\kappa} 
\|\bar G\|_{L_{p}(C_{r _{1}})}.
\medskip$$

\noindent
It follows by Lemma \ref{lemma 10.23.1} that  
$$
[v]_{C^{\kappa }(  C_{r'_{0}}(t_{0},x_{0}))}\leq N_{1}
I(\theta_{0},\varepsilon,r_{1}),
$$
where   $N_{1}$ depends only on $d$  and $\delta$. 
We can now specify 
$\theta_{0}$ and $\varepsilon$. First we chose
$\varepsilon\in(0,1)$ so that
$$
N_{1}\varepsilon^{ \kappa_{0}-\kappa }
(1-\varepsilon)^{-\kappa_{0}}= 1/4.
\smallskip$$
Since $\kappa_{0}-\kappa\geq (\kappa_{0}-1)/2>0$ and $\kappa_{0}$  
depends only on $d$ and $\delta$ and $N_{1}$
depends only on $d$  and $\delta$,
 $\varepsilon$ also 
depends only on $d$  and $\delta$.
After that we take   $\theta_{0}=\theta_{0}(d, 
\delta)\in(0,1]$ so that $N_{1}\theta_{0}
\varepsilon ^{-\kappa }\leq1/4$.

Then
\begin{equation}
                                              \label{010.25.1}
[v]_{C^{\kappa }(  C_{r'_{0}}(t_{0},x_{0}))}\leq 
(1/2)[v]_{C^{  \kappa}(  C_{r_{1}})}
+NJ,
\medskip\end{equation}
where $N=N(d,\delta)$ and   \smallskip
$$
J=K_{0}\sup_{ C_{r_{1}}}\big(|v|+|Dv|\big)
+  K_{1}+\|\bar G\|_{L_{p}(C_{r _{1}})}.
$$

Now observe that if $(t,x),(s,x)\in C_{r_{0}}$ and $t> s$, then
 either $|t-s|<\big(r_{0}'\big)^{2}$, in which case 
$(t,x)\in C_{r'_{0}}(s,x)$ and \vspace{5pt}
$$
(t-s)^{-\kappa /2}\big|v(t,x)-v(s,x)\big|
\leq (1/2)[v]_{C^{  \kappa}(  C_{r_{1}})}
+NJ
\medskip$$
owing to \eqref{010.25.1},
 or
  $|t-s|\geq \big(r_{0}'\big)^{2}$ when   \vspace{5pt}
\begin{align*}
\big|v(t,x)-v(s,x)\big|\leq &\, 2(t-s)^{\kappa /2}(r'_{0})^{- \kappa }
\sup_{  C_{r_{1}}}|v|
\\[10pt]
\leq  &\, N(t-s)^{\kappa /2}(r_{1}-r_{0})
^{- \kappa }\sup_{  C_{r_{1}}}|v|.
\end{align*}

Next if $(t,x),(t,y)\in C_{r_{0}}$ and $x\ne y$, then either
$|x-y|<r_{0}'$, in which case $(t,y)\in  C_{r'_{0}}(t,x)$ and  \smallskip
$$
|x-y|^{-(\kappa-1) }\big|Dv(t,x)-Dv(t,y)\big|\leq
(1/2)[v]_{C^{  \kappa}(  C_{r_{1}})}+NJ,
$$
or   $|x-y|\geq r_{0}'$ and  \smallskip
\begin{align*}
|Dv(t,x)-Dv(t,y)|\leq &\, 2|x-y|^{\kappa-1 }\big(r_{0}'\big)^{-(\kappa-1) }
\sup_{  C_{r_{1}}}|Dv|
\\[10pt]
\leq &\, N|x-y|^{\kappa-1 }
(r_{1}-r_{0})
^{-(\kappa-1) }\sup_{  C_{r_{1}}}|Dv|.
\end{align*}
This proves \eqref{1,3.5} and the lemma.\qed

\begin{theorem}
                                  \label{theorem 1,3.4}
Take $0<r<R\leq R_{0}$ and take  
$\kappa $  and $\theta_{0}$
from Lemma \ref{lemma 1,3.4}
and suppose that Assumption \ref{assumption 12,29.2} \(iv\,\)
 is satisfied with this
$\theta_{0}$.
Let 
$ v \in W^{1,2}_{\infty}(C_{R}) $ be a  
  solution
of \eqref{9.28.4} in $C_{R}$
 such that
$$
\osc_{  C_{R}} v< \hat R_{0},
$$

 Then 
\begin{equation}
                                                \label{3,14,3}
[v]_{C^{\kappa}(  C_{r})}\leq 
N(R-r)^{-\kappa}\sup_{  C_{R}}|v|
+N\big( K_{1}+\|\bar G\|_{L_{p}(C_{R})}\big) ,
\medskip\end{equation}
where $N$ depends only on $d,\delta$, and $K_{0}$.

\end{theorem}

Proof. Fix a number $c\in(0,1)$ such that $c^{4}>3/4$ and   introduce
$$
r_{0}=r,\quad
r_{n}=r+c_{0}(R-r)\sum_{k=1}^{n}c^{k},\quad n\geq1,
$$
where $c_{0}$ is chosen in such a way that $r_{n}\to R$ as $n\to\infty$.
Then 
Lemma \ref{lemma 1,3.4} and \eqref{11.1.1} allow us to find constants
$N_{1}$ and $N$ depending only on $d,\delta$, and $K_{0}$,
  such that for all $n$
and $\varepsilon\in(0,1)$   \smallskip
\begin{align*}
A_{n}:=[v]_{C^{\kappa}(  C_{r_{n}})}\leq &\, 
\big(2^{-1}+N_{1}\varepsilon^{\kappa-1}\big)A_{n+2}
+N( K_{1}+\|\bar G\|_{L_{p}(C_{R})})
\\[10pt]
&\,  +
N(R-r)^{-\kappa}
c^{-n\kappa}(1+\varepsilon^{-1})\sup_{  C_{R}}|v|.
\end{align*}
We choose $\varepsilon< 1$ so  
that $2^{-1}+N_{1}\varepsilon^{\kappa-1}
\leq3/4$ and then recalling that $\kappa\leq 2$ conclude that  \smallskip
\begin{align*}
\sum_{m=0}^{\infty}(3/4)^{m}A_{2m}\leq
\sum_{k=1}^{\infty}(3/4)^{m}A_{2m} +N( K_{1}+\|\bar G\|_{L_{p}(C_{R})})
\\[10pt]
+N(R-r)^{-\kappa}\sup_{  C_{R}}|v|
\sum_{m=0}^{\infty}(3/4)^{m} c^{-4m},
\end{align*}
where the last series converges since $3c^{-4}/4<1$.
By canceling like terms we come to \eqref{3,14,3}
and the theorem is proved. \qed

{\bf Proof of Theorem          
\ref{theorem 12,30.1}}.  
 We keep assuming that
 that $g\in W^{1,2}_{\infty}(\bR^{d+1})$ and
$\bar G$ is bounded. 

Since
$$
|H(\sfu',0,t,x)|=|G(\sfu',0,t,x)|\leq K_{0}|\sfu'|+\bar G(t,x),
$$
  for $H_{K}=\max(H,P-K)$ we have 
$$
|H_{K}(\sfu',0,t,x)|\leq|H(\sfu',0,t,x)|\leq K_{0}|\sfu'|+\bar G(t,x).
$$
It follows by Lemma 3.2 of \cite{Kr_17.1} that for any $K$
there exist measurable $\bS_{\bar\delta}$-valued $a$, $\bR^{d}$-valued
$b$, and real-valued $f$ such that in $\Pi$ (a.s.)
\begin{equation}
                                                 \label{10.1.1}
\partial_{t}v_{K}+a^{ij}D_{ij}v_{K}+b^{i}D_{i}v_{K}+f=0
\end{equation}
and $|b|\leq K_{0}$, $|f|\leq \bar G+K_{0}|v_{K}|$.

Then by the parabolic Aleksandrov estimates   \smallskip   
\begin{equation}
                                                  \label{3,14,1} 
|v_{K}|\leq N (\|g\|_{C (\Pi)}+\|\bar{G}\|_{L_{d+1}(\Pi)}
  ),
\medskip\end{equation}
where $N$ depends only on  
$d$, $\delta$, $K_{0}$, $T$, and the diameter
of $\Omega$. 

Then, since \eqref{10.1.1} has form of a linear equation, by the
well-known results from the linear theory
we  estimate  not only 
$|v_{K}|$  but also  the modulus of continuity
of $v_{K}$ through that of $g$, $\sup|g|$, and $\|\bar G\|_{L_{d+1}(\Pi)}$
with constants independent of $K$.

Hence, the family $\{v_{K};K\geq1\}$
is equicontinuous on $\bar\Pi$. More precisely
there exists a function $\bar\omega(\varepsilon)$,
depending only on $\varepsilon$, 
 $
 \|\bar G\|_{L_{d+1}(\Pi )}  
 $,
 $\delta$, $d$, $K_{0}$, $\rho_{\ext}(\Omega)$, $ \|g\|_{C (\Pi)}$,
  and the modulus   of
continuity of $g$ on $\partial'\Pi $,
such that $\bar\omega(\varepsilon)\to0$ as 
$\varepsilon\downarrow 0$ and  \smallskip
\begin{equation}
                                                 \label{3,14,2} 
\big|v_{K}(t,x)-v_{K}(s,y)\big|\leq\bar\omega\Big(\rho\big((t,x),(s,y)\big)
\Big)
\medskip\end{equation}
for any $x,y\in\bar\Omega$ and $s,t\in[0,T]$.

It follows that there is a sequence $K_{n}\to\infty$ and a function
$v$ such that $v^{n}:=v_{K_{n}}\to v$ uniformly in $\bar\Pi$.  
Of course, \eqref{1.25.6} holds,
owing to Theorem \ref{theorem 1,3.4},
which also implies that   $Dv^{n}\to Dv$         \label{page 3,22,1}   
locally uniformly in $\Pi$.

The following lemma, in which the boundedness of $\bar G$ is not used,   
 will allow us to prove that
$v$ is a  viscosity  solution. Introduce  \smallskip
$$
H^{0}(\sfu'',t,x)=H\big(v(t,x),Dv(t,x),\sfu'',t,x\big).
\smallskip$$
\begin{lemma}
                                           \label{lemma 3,14,3}  
There is a constant $N$, depending only on $d$ and
$\delta$, such that
for any $C_{r}(t,x)$ satisfying $ \bar C_{r}(t,x)\subset
\Pi$  and
$\phi\in W^{1,2}_{d+1}(C_{r}(t,x))$ we have on $C_{r}(t,x)$ that \smallskip
\begin{equation}
                                             \label{3,14,5}
v\leq \phi+Nr^{d/(d+1)}\big\|\big(\partial_{t}\phi+
H^{0}[\phi]\big)_{+}\big\|_{L_{d+1}(C_{r}(t,x))}
+\max_{\partial'C_{r}(t,x)}(v-\phi)_{+},
\medskip\end{equation}
\begin{equation}
                                             \label{3,14,6}
v\geq \phi-Nr^{d/(d+1)}\big\|\big(\partial_{t}\phi+
H^{0}[\phi] \big)_{-}\big\|_{L_{d+1}(C_{r}(t,x))}
-\max_{\partial'C_{r}(t,x)}(v-\phi)_{-} .
\end{equation}
\end{lemma}

Proof. For $m=1,2,...$ introduce  \smallskip
$$
H^{m}(\sfu'',t,x)=\sup_{n\geq m}
H\big(v^{n}(t,x),Dv^{n}(t,x),\sfu'',t,x\big)
$$
and observe that for $n\geq m$  \smallskip
$$
\partial_{t}v_{n}+\max\big(H^{m}[v^{n}], P[v^{n}]-K_{n}\big)\geq0,
$$
implying that  \smallskip
\begin{align*}
-\partial_{t}\phi-\max\big(H^{m}[\phi] ,P[\phi]-K_{n}\big)\leq &\,
-\partial_{t}\phi-\max\big(H^{m}[\phi] ,P[\phi]-K_{n}\big)
\\[10pt]
&\, +\partial_{t}v^{n}+
\max\big(H^{m}[v^{n}] ,P[v^{n}]-K_{n}\big) 
\\[10pt]
= &\,
\partial_{t}(v^{n}-\phi)+a^{ij}D_{ij}(v^{n}-\phi),
\end{align*}
where $a=(a^{ij})$ is an $\bS_{\bar{\delta}}$-valued
function.

It follows by the parabolic Aleksandrov estimates 
that     \smallskip
$$
v^{n}\leq \phi+\max_{\partial'C_{r}(t,x)}(v^{n}-\phi)_{+}
\medskip$$
$$
+Nr^{d/(d+1)}\big\|\big\{\partial_{t}\phi+ 
\max\big(H^{m}[\phi] ,P[\phi]-K_{n}\big)\big\}_{+}\big\|_{L_{d+1}(C_{r}(t,x))},
\medskip$$
where   $N=N(d,\delta)$.  By sending $n\to\infty$
and using the dominated convergence theorem,
we obviously get    \smallskip
\begin{equation}
                                             \label{3,14,4}
v \leq \phi+\max_{\partial'C_{r}(t,x)}(v -\phi)_{+}
+Nr^{d/(d+1)}\big\|\big(\partial_{t}\phi+ 
  H^{m}[\phi]  \big)_{+}\big\|_{L_{d+1}(C_{r}(t,x))}.
\medskip\end{equation}

Since $v^{n}\to v$ and $Dv^{n}\to Dv$ uniformly
in $C_{r}(t,x))$ we have that $\big|H^{m}[\phi] \big|\leq 
\big|H^{m}[0] \big|+N|D^{2}\phi|$, where $\big|H^{m}[0] \big|$
is dominated  by an $L_{d+1}$-function independent of $m$.   
Furthermore, the continuity of $H(\sfu,t,x)$ 
with respect to $\sfu$ implies that
$H^{m}[\phi] \to H^{0}[\phi] $ as $m\to\infty$
at any point in $C_{r}(t,x) $. By using
the dominated convergence theorem once more
and sending $m\to\infty$ in \eqref{3,14,4}
we come to \eqref{3,14,5}.

Similarly \eqref{3,14,6} is established.
The lemma is proved. \qed

Now we prove that $v$ is  a   $C$-viscosity solution   
of \eqref{7.29.10} if part (a) of Assumption \ref{assumption 12,29.1} (ii)
is satisfied.
 Let $(t_{0},x_{0})\in\Pi$, $r>0$, and 
$\phi\in C^{1,2} \big(C_{r}(t_{0},x_{0})\big)$
 be such that $\bar C_{r}(t_{0},x_{0})\subset \Pi$
and $v-\phi$ attains a local
maximum at $(t_{0},x_{0})$.
 Then for   $\varepsilon
>0$ and all small $r>0$ for   \smallskip
$$
\phi_{\varepsilon,r}(t,x)=\phi (t,x)-\phi(t_{0},x_{0})
+v(t_{0},x_{0})+\varepsilon(
|x-x_{0}|^{2}+t-t_{0}- r^{2})
$$
 we have that 
$$
\max_{\partial'C_{r}(t_{0},x_{0})}(v -\phi_{\varepsilon,r})_{+}
=0.
$$
Hence, by Lemma \ref{lemma 3,14,3}
\begin{align*}
   \varepsilon r^{2}= 
(v -\phi_{\varepsilon,r})(t_{0},x_{0})
\leq Nr^{d/(d+1)}\big\|\big(\partial_{t}\phi_{\varepsilon,r}+
H^{0}[\phi_{\varepsilon,r}]\big)_{+}\big\|_{L_{d+1}(C_{r}(t_{0},x_{0}))}
\\[10pt]
=Nr^{d/(d+1)}\big\|\big(\partial_{t}\phi_{\varepsilon }+
H^{0}[\phi_{\varepsilon }]\big)_{+}\big\|_{L_{d+1}(C_{r}(t_{0},x_{0}))},
\end{align*}
where $\phi_{\varepsilon }=\phi+\varepsilon
\big( |x |^{2}+ t\big)$. It follows that  \smallskip
$$
Nr^{-(d+2)}\big\|\big(\partial_{t}\phi_{\varepsilon }+
H^{0}[\phi_{\varepsilon }]\big)_{+}\big\|^{d+1}_{L_{d+1}(C_{r}(t_{0},x_{0}))}
\geq \varepsilon^{d+1},
\medskip$$
$$
N\esssup_{C_{r}(t_{0},x_{0}) }
\big(\partial_{t}\phi_{\varepsilon }+
H^{0}[\phi_{\varepsilon }]\big)\geq \varepsilon,
\medskip$$
$$
N\lim _{ r\downarrow0}\esssup_{C_{r}(t_{0},x_{0}) }
\big(\partial_{t}\phi_{\varepsilon }+
H^{0}[\phi_{\varepsilon }]\big)\geq \varepsilon.
\medskip$$
By letting $\varepsilon\downarrow0$ we conclude that  \medskip
\begin{equation}
                                                 \label{3.22.1}
 \lim _{ r\downarrow0}\esssup_{C_{r}(t_{0},x_{0}) }
\big[\partial_{t}\phi(t,x) +
H\big(v(t,x),Dv(t,x),D^{2}\phi(t,x),t,x\big)\big]\geq 0.
\end{equation} 

Now, note that, as $(t,x)\to(t_{0},x_{0})$, we have
$v(t,x)\to v(t_{0},x_{0})$ and (see Remark \ref{remark 12,14.2})  
$Dv(t,x)\to Dv(t_{0},x_{0})$. Also 
 $
D\phi(t,x)\to D\phi(t_{0},x_{0})=Dv(t_{0},x_{0})
 $ and $
D^{2}\phi(t,x)\to D^{2}\phi(t_{0},x_{0}) 
 $.
It follows by Assumption \ref{assumption 12,29.1} (ii) (a)
that one can replace $Dv(t,x)$ in \eqref{3.22.1} with $D\phi(t,x)$.
Then, so modified 
 \eqref{3.22.1} implies \eqref{1,29,2} meaning that $v$   
is a  $C$-viscosity subsolution.

The fact that it is also a   $C$-viscosity supersolution  
is proved similarly on the basis of \eqref{3,14,6}.

In case Assumption \ref{assumption 12,29.1} (ii) (b)
is satisfied, we still come to \eqref{3.22.1}
and can replace $Dv(t,x)$   with $D\phi(t,x)$
just because of the continuity of $H(\sfu,t,x)$ in $[u']$
uniform with respect to $\sfu'_{0}$, $\sfu''$, and $(t,x)$.

This proves Theorem \ref{theorem 12,30.1}
in our particular case that  
$g\in W^{1,2}_{\infty}(\bR^{d+1})$ and $\bar G$ is bounded.
In the same case we also have estimates
\eqref{3,14,1}, \eqref{3,14,2} with $v$ in place of $v_{K}$
and the same $N$ and $\bar\omega$. Also \eqref{1.25.6} holds
by the above.

In the case of general $\bar G$, for $n=1,2,...$,
we replace $H(\sfu,t,x)$ in \eqref{7.29.10} with  \smallskip
$$
H(\sfu,t,x)I_{\bar{G}(t,x)\leq n}
+F(\sfu'_{0},\sfu'',t,x)I_{\bar{G}(t,x)> n}
$$ \smallskip
$$
=F(\sfu'_{0},\sfu'',t,x)+G(\sfu,t,x)I_{\bar{G}(t,x)\leq n}
\medskip$$
and apply the already proved  version of Theorem \ref{theorem 12,30.1}
to introduce $u_{n}$ as the
$C$-viscosity or $L_{p}$-viscosity solutions
  in $\Pi$ of so obtained equations 
 with the same boundary
condition $u_{n}=g$ on $\partial'\Pi$.

From the above we see that the estimates 
\eqref{3,14,1} and \eqref{3,14,2} 
  with $u_{n}$ in place of $v_{K}$  
 and the estimates    
of   $[u_{n}]_{C^{\kappa}(  C_{r}(t,x))}$
 are uniform with
respect to $n$. This and the fact that
the boundedness of $\bar G$ is not used
in Lemma \ref{lemma 3,14,3}  allow  us to repeat  what was said  
about $v^{n}$ with obvious changes and 
proves the theorem for general $H$ but 
still assuming that $g\in W^{1,2}_{\infty}(\bR^{d+1}) $.

One drops this assumption by using uniform
approximations of $g$ by smooth ones
preserving the modulus of continuity on
$\partial'\Pi$. This guarantees that
for the approximating solutions the estimates like 
\eqref{3,14,1} and \eqref{3,14,2}
will hold and referring to the argument in the previous
paragraph
brings the proof
of Theorem \ref{theorem 12,30.1} to an end.   \qed


\begin{thebibliography}{mm}

\bibitem{Caf89} L.A. Caffarelli, {\em Interior a priori estimates 
for solutions of fully non-linear equations\/},  Ann. Math.,  
Vol. 130 (1989), 189--213.

\bibitem{CIL} M.G. Crandall, H. Ishii, and P.-L. Lions,
{\em
User's guide to viscosity solutions of second order 
partial differential equations\/},
Bull. Amer. Math. Soc. (N.S.), Vol. 27 (1992), 1--67.

\bibitem{CKS00} M. G. Crandall, M. Kocan, and A. \'Swi{\c e}ch, {\em
$L^p$-theory for fully nonlinear uniformly parabolic equations\/},
 Comm. Partial Differential Equations, Vol. 25  (2000),
 No. 11-12, 1997--2053.

\bibitem{ST_16} J.V. da Silva and E.V. Teixeira,
{\em
Sharp regularity estimates for second order fully
nonlinear parabolic equations\/}, Math. Ann.
DOI 10.1007/s00208-016-1506-y

\bibitem{Go_62} K. K. Golovkin, 
{\em On equivalent normalizations of fractional spaces\/},
Automatic programming, numerical methods and functional analysis, Trudy
Mat. Inst. Steklov., Vol. 66, Acad. Sci. USSR, Moscow--Leningrad, 1962, 
364--383.


\bibitem{Kr_96}  N.V. Krylov, ``Lectures on  elliptic and parabolic
equations in H\"older spaces'', Amer.
Math. Soc., Providence, RI, 1996; Russian translation,
``Nauchnaya kniga'', Novosibirsk, 1998.


\bibitem{Kr_13.1} N.V. Krylov, {\em On the existence of $W^{2}_{p}$ 
solutions for fully nonlinear elliptic
equations under  relaxed convexity assumptions\/}, 
Comm. Partial Differential Equations, Vol. 38 (2013), No. 4,  
687--710.

\bibitem{Kr_17.1} N.V. Krylov, {\em
On the existence of $W^{1,2}_{p}$ 
solutions for fully nonlinear parabolic
equations under either relaxed 
or no  convexity assumptions\/},  
Harvard University, Center of Mathematical Sciences and 
Applications,
 Nonlinear Equation
Publication, http://arxiv.org/abs/1705.02400

 \bibitem{LSU} O.A. Ladyzhenskaya, V.A. Solonnikov, and
N.N. Ural'tseva,
 ``Linear and quasi-linear parabolic equations'',
Nauka, Moscow, 1967, in Russian; English translation: 
Amer. Math. Soc.,
Providence, RI, 1968.

\bibitem{Li} G.M. Lieberman,
 ``Second order parabolic
 differential equations'', World
Scientific, Singapore, 1996.

\bibitem{Sa_88}  M. V. Safonov, {\em
On the classical solutions of nonlinear elliptic
equations of second order\/}, Izvestija Acad. Nauk SSSR, ser. matemat.,
Vol. 52 (1988), No. 6, 1272--1287   in Russian; English translation
in Math. USSR Izvestiya, Vol. 33 (1989), No. 3, 597--612.
 

\bibitem{SS_2014} L. Silvestre and B. Sirakov,
{\em Boundary regularity for viscosity
solutions of fully nonlinear elliptic equations\/}, Comm. Partial
Differential Equations, Vol. 39 (2014), No. 9, 1694--1717.

\bibitem{ST_2015} L. Silvestre and E.V. Teixeira, 
{\em Regularity estimates for fully non
linear elliptic equations which are asymptotically convex\/},
 Contributions
to nonlinear elliptic equations and systems, 425--438, Progr. Nonlinear
Differential Equations Appl., Vol. 86, Birkh\"auser/Springer, Cham, 2015.

\bibitem{Sw_97} A. \'Swi{\c e}ch,
{\em $W^{1}_{ p}$-interior estimates for solutions of fully 
nonlinear, uniformly elliptic equations\/},
 Adv. Differential Equations, Vol. 2 (1997), 1005--1027.

\bibitem{Tr_88} N.S.
Trudinger, {\em Comparison principles and pointwise estimates
for viscosity solutions of nonlinear elliptic equations\/},
 Rev.
Mat. Iberoamericana, Vol. 4 (1988), No. 3--4, 453--468.

\bibitem{Tr_89} N.S.
Trudinger, {\em On regularity and existence of viscosity
solutions of nonlinear second order, elliptic equations\/},
 Partial
differential equations and the calculus of variations, Vol. II,
939--957, Progr. Nonlinear Differential Equations Appl., Vol. 2,
Birkh\"auser Boston, Boston, MA, 1989. 

\bibitem{Wa92_1}  L. Wang, {\em
On the regularity of fully nonlinear
parabolic equations: II\/},  Comm. Pure Appl. Math., Vol. 45 (1992),
 141--178.


\end{thebibliography}
\end{document}